\pdfoutput=1
\newif\ifpersonal

\personaltrue 

\ifpersonal
\newcommand*{\personal}[1]{\tcolor[rgb]{0,0,1}{(Personal: #1)}}
\newcommand*{\todo}[1]{\textcolor{red}{(Todo: #1)}}
\else
\newcommand*{\personal}[1]{\ignorespaces}
\newcommand*{\todo}[1]{\ignorespaces}
\fi

\documentclass[oneside, hidelinks, 12pt, a4paper,reqno]{amsart}
\usepackage{extsizes}
\usepackage[toc,page]{appendix}
\usepackage[foot]{amsaddr}
\usepackage{etoolbox}

\makeatletter
\patchcmd{\@setaddresses}{\indent}{\noindent}{}{}
\patchcmd{\@setaddresses}{\indent}{\noindent}{}{}
\patchcmd{\@setaddresses}{\indent}{\noindent}{}{}
\patchcmd{\@setaddresses}{\indent}{\noindent}{}{}
\makeatother
\usepackage{amsthm,mathtools,bm,tensor}

\usepackage[T1]{fontenc}
\DeclareFontFamily{U}{min}{}
\DeclareFontShape{U}{min}{m}{n}{<-> dmjhira}{}
\usepackage[bitstream-charter]{mathdesign}
\usepackage{XCharter}
\let\mathcal\undefined
\DeclareMathAlphabet{\mathcal}{U}{dutchcal}{m}{n}
\usepackage{enumerate,comment,braket,xspace,tikz-cd,csquotes}
\usepackage[driver=dvips,centering,vscale=0.7,hscale=0.8]{geometry}
\usepackage{url}
\usepackage{microtype}
\usepackage[greek.polutoniko,english]{babel}
\usepackage[style=authoryear,backend=bibtex8,sorting=nyt,citestyle=ieee-alphabetic, bibstyle=alphabetic,texencoding=ascii,bibencoding=utf8,maxnames=5,minnames=3,maxbibnames=99]{biblatex}
\usepackage{setspace}
\usepackage{enumitem}
\usepackage[cal=dutchcal, scr=boondoxo]{mathalfa}
\usepackage{listings}
\usepackage{colordvi}
\usepackage{tikz,pgf}
\usepackage[normalem]{ulem}
\usetikzlibrary{arrows,decorations.pathreplacing,backgrounds,positioning,fit,matrix,calc}
\tikzset{curve/.style={settings={#1},to path={(\tikztostart)
			.. controls ($(\tikztostart)!\pv{pos}!(\tikztotarget)!\pv{height}!270:(\tikztotarget)$)
			and ($(\tikztostart)!1-\pv{pos}!(\tikztotarget)!\pv{height}!270:(\tikztotarget)$)
			.. (\tikztotarget)\tikztonodes}},
	settings/.code={\tikzset{quiver/.cd,#1}
		\def\pv##1{\pgfkeysvalueof{/tikz/quiver/##1}}},
	quiver/.cd,pos/.initial=0.35,height/.initial=0}
\tikzset{tail reversed/.code={\pgfsetarrowsstart{tikzcd to}}}
\tikzset{2tail/.code={\pgfsetarrowsstart{Implies[reversed]}}}
\tikzset{2tail reversed/.code={\pgfsetarrowsstart{Implies}}}
\usepackage{epsfig}
\usepackage{epstopdf}
\usepackage{tikz-cd}
\usepackage{fancyhdr}
\usepackage{graphicx}
\usepackage{fullpage}
\makeatletter
\newcommand{\mylabel}[2]{#2\def\@currentlabel{#2}\label{#1}}
\makeatother
\usepackage{hyperref}
\hypersetup{
	colorlinks=true,
	linkcolor=red,
        linktoc=page,
	anchorcolor=black,
	filecolor=blue,
	citecolor = blue, 
	menucolor=., 
	urlcolor=cyan,
	breaklinks=true,
    unicode=true,
}
\usepackage[capitalize]{cleveref}
\usepackage{nameref}
\usepackage{enumitem}
\crefalias{enumi}{lemman}
\usepackage{comment}
\usepackage{color}
\usepackage[totoc]{idxlayout}
\usepackage[lastexercise]{exercise}
\usepackage{xcolor}
\usepackage{breakcites}

\makeatletter
\def\@tocline#1#2#3#4#5#6#7{\relax
	\ifnum #1>\c@tocdepth 
	\else
	\par \addpenalty\@secpenalty\addvspace{#2}%
	\begingroup \hyphenpenalty\@M
	\@ifempty{#4}{%
		\@tempdima\csname r@tocindent\number#1\endcsname\relax
	}{%
		\@tempdima#4\relax
	}%
	\parindent\z@ \leftskip#3\relax \advance\leftskip\@tempdima\relax
	\rightskip\@pnumwidth plus4em \parfillskip-\@pnumwidth
	#5\leavevmode\hskip-\@tempdima
	\ifcase #1
	\or\or \hskip 1em \or \hskip 2em \else \hskip 3em \fi%
	#6\nobreak\relax
	\dotfill\hbox to\@pnumwidth{\@tocpagenum{#7}}\par
	\nobreak
	\endgroup
	\fi}
\makeatother

\makeatletter
\renewcommand\paragraph{\@startsection{paragraph}{5}{\z@}%
  {3.25ex \@plus1ex \@minus.2ex}%
  {-1em}%
  {\normalfont\normalsize\bfseries}}
\renewcommand\subparagraph{\@startsection{subparagraph}{6}{\parindent}%
  {3.25ex \@plus1ex \@minus .2ex}%
  {-1em}%
  {\normalfont\normalsize\bfseries}}
\def\toclevel@subsubsubsection{4}
\def\toclevel@paragraph{5}
\def\toclevel@paragraph{6}
\def\l@subsubsubsection{\@dottedtocline{4}{7em}{4em}}
\def\l@paragraph{\@dottedtocline{5}{10em}{5em}}
\def\l@subparagraph{\@dottedtocline{6}{14em}{6em}}
\makeatother

\usepackage{marginnote}

\usepackage{lipsum}
\usepackage{fancyhdr}
\providecommand{\abstract}{}
\usepackage{abstract}
\DeclareGraphicsExtensions{.png,.jpg,.pdf,.mps}
\usepackage{fancyhdr}
\usepackage{lastpage}
\usepackage{multirow}
\usepackage{faktor}

\usepackage{array}
\usepackage{makecell}

\newcommand{\stackspace}{1.7}
\newcommand{\stack}[2][1cm]{\;\tikz[baseline, yshift=.65ex]%
	{\foreach \k [evaluate=\k as \r using (.5*#2+.5-\k)*\stackspace] in {1,...,#2}{%
			\ifodd\k{\draw[->](0,\r pt)--(#1,\r pt);}%
			\else{\draw[<-](0,\r pt)--(#1,\r pt);}\fi
	}}\;}

\tikzcdset{
  diagrams={>={Straight Barb[scale=0.5]}}
}
\tikzset{
  altstackar/.style={decorate, decoration={show path construction,
    lineto code={
      \path (\tikzinputsegmentfirst); \pgfgetlastxy{\xstart}{\ystart}
      \path (\tikzinputsegmentlast); \pgfgetlastxy{\xend}{\yend}
      \path ($(0,0)!1.5pt!(\ystart-\yend,\xend-\xstart)$); \pgfgetlastxy{\xperp}{\yperp}
      \foreach \n[evaluate=\n as \k using .5*#1-\n+.5] in {1,...,#1}{
        \ifodd\n{\draw[->, shorten <=2pt, shift={($\k*(\xperp,\yperp)$)}](\xstart,\ystart)--(\xend,\yend);}
        \else{\draw[<-, shorten >=2pt, shift={($\k*(\xperp,\yperp)$)}](\xstart,\ystart)--(\xend,\yend);}\fi
      }
    }
  }}, altstackar/.default={1}
}

\usepackage{pict2e}

\makeatletter
\newcommand{\adjunction}[4]{%
  #1\colon #2%
  \mathrel{\vcenter{%
    \offinterlineskip\m@th
    \ialign{%
      \hfil$##$\hfil\cr
      \longrightharpoonup\cr
      \noalign{\kern-.3ex}
      \smallbot\cr
      \longleftharpoondown\cr
    }%
  }}%
  #3 \noloc #4%
}
\newcommand{\longrightharpoonup}{\relbar\joinrel\rightharpoonup}
\newcommand{\longleftharpoondown}{\leftharpoondown\joinrel\relbar}
\newcommand\noloc{%
  \nobreak
  \mspace{6mu plus 1mu}
  {:}
  \nonscript\mkern-\thinmuskip
  \mathpunct{}
  \mspace{2mu}
}
\newcommand{\smallbot}{%
  \begingroup\setlength\unitlength{.15em}%
  \begin{picture}(1,1)
  \roundcap
  \polyline(0,0)(1,0)
  \polyline(0.5,0)(0.5,1)
  \end{picture}%
  \endgroup
}
\makeatother

\usepackage[capitalize]{cleveref}
\usepackage{aliascnt}

\theoremstyle{plain}
\newtheorem{theorem}{Theorem}[section]

\newaliascnt{conjecture}{theorem}

\aliascntresetthe{conjecture}

\theoremstyle{definition}
\newaliascnt{definition}{theorem}
\newtheorem{definition}[definition]{Definition}
\aliascntresetthe{definition}

\newaliascnt{example}{theorem}
\newtheorem{example}[example]{Example}
\aliascntresetthe{example}

\newaliascnt{remark}{theorem}
\newtheorem{remark}[remark]{Remark}
\aliascntresetthe{remark}

\newaliascnt{construction}{theorem}
\newtheorem{construction}[construction]{Construction}
\aliascntresetthe{construction}

\theoremstyle{plain}
\newaliascnt{lemma}{theorem}
\newtheorem{lemma}[lemma]{Lemma}
\aliascntresetthe{lemma}

\newaliascnt{corollary}{theorem}
\newtheorem{corollary}[corollary]{Corollary}
\aliascntresetthe{corollary}

\newaliascnt{proposition}{theorem}
\newtheorem{proposition}[proposition]{Proposition}
\aliascntresetthe{proposition}

\newaliascnt{assumptions}{theorem}

\aliascntresetthe{assumptions}

\Crefname{theorem}{Theorem}{Theorems}
\Crefname{conjecture}{Conjecture}{Conjectures}
\Crefname{definition}{Definition}{Definitions}
\Crefname{example}{Example}{Examples}
\Crefname{remark}{Remark}{Remarks}
\Crefname{construction}{Construction}{Constructions}
\Crefname{lemma}{Lemma}{Lemmas}
\Crefname{corollary}{Corollary}{Corollaries}
\Crefname{proposition}{Proposition}{Propositions}
\Crefname{assumptions}{Assumption}{Assumptions}

 \newcommand{\cB}{\mathscr B}
 \newcommand{\cC}{\mathscr C}

 \newcommand{\cF}{\mathscr F}

 \newcommand{\cV}{\mathscr V}

\newcommand{\CAlg}{\mathrm{CAlg}}

\newcommand{\Dafd}{\mathrm{dAfd}}
\newcommand{\Daff}{\mathrm{dAff}}
\newcommand{\Afd}{\mathrm{Afd}}
\newcommand{\Aff}{\mathrm{Aff}}
\newcommand{\Sp}{\mathrm{Sp}}
\newcommand{\Spec}{\mathrm{Spec}}
\newcommand{\Spf}{\mathrm{Spf}}

\newcommand{\rig}{\mathrm{rig}}

\newcommand{\Perf}{\mathrm{Perf}}
\newcommand{\cat}{$\infty$-category }
\newcommand{\cates}{$\infty$-categories }
\newcommand{\catperf}{\mathrm{Cat}^{\Perf}}

\newcommand{\Cohb}{\mathrm{Coh}^\mathrm{b}}

\newcommand{\Fun}{\mathrm{Fun}}
\newcommand{\Hom}{\mathrm{Hom}}
\newcommand{\Mod}{\mathrm{Mod}}
\newcommand{\LMod}{\mathrm{LMod}}
\newcommand{\Ind}{\mathrm{Ind}}

\newcommand{\cofib}{\mathrm{cofib}}

\newcommand{\m}{\mathscr}
\newcommand{\f}{\mathrm }
\newcommand{\F}{\mathscr}

\newcommand{\pr}{\mathrm{Pr}^{\mathrm{L}}}
\newcommand{\prst}{\mathrm{Pr}^{\mathrm{L}}_{\mathrm{st}}}
\newcommand{\promega}{\mathrm{Pr}^{\mathrm{L}, \omega}}
\newcommand{\promegast}{\mathrm{Pr}^{\mathrm{L},\omega}_{\mathrm{st}}}
\newcommand{\Addresses}{{
  \bigskip
  \footnotesize

  Matteo Montagnani, \textsc{SISSA, Via Bonomea 265, 34136 Trieste TS, Italy}\par\nopagebreak
  \textit{E-mail address}: \href{mailto:mmontagn@sissa.it}{\texttt{mmontagn@sissa.it}}

  \medskip

}}

\usepackage{doc}

\title{An algebraicity criterion for rigid analytic varieties and formal schemes via perfect complexes}
\author{Matteo Montagnani}
\begin{document}

\maketitle
\section*{Abstract}
In this paper, we extend a theorem of To\"en and Vaqui\'e to the non-Archimedean and formal settings. More precisely, we prove that a smooth and proper rigid analytic variety is algebraizable if and only if its category of perfect complexes is smooth and proper. As a corollary, we deduce an analogous statement for formal schemes and demonstrate that, in general, the bounded derived category of coherent sheaves on a formal scheme is not smooth.
\tableofcontents 

\section*{Introduction}
In algebraic geometry, it is standard practice to study an algebraic variety $X$ through its abelian category of quasi-coherent sheaves, as well as its derived category, commonly referred to as the derived triangulated category of the algebraic variety (see \cite{bondal1,orlov}). However, for our purposes, it is more natural to focus on an enhanced version of this category, which we denote by $\mathrm{QCoh}(X)$. This enhanced category has been extensively studied in the context of dg-categories, for example in \cite{toen2007homotopy, canonaco2017tour}, and more recently within the framework of \cates (see \cite{SAG,antieau2018uniquenesantis}).

The derived category and its enhancements are important and well-studied geometric invariants. In \cite{bondal2002generators}, it was proven that, under mild conditions on a scheme $X$, the derived category of quasi-coherent sheaves admits a single compact generator. Moreover, in the $\infty$-categorical setting, if $X$ is an Artin stack, then $\mathrm{QCoh}(X)$ admits both a t-structure and a symmetric monoidal structure. Using these structures, if $X$ is a Deligne-Mumford (DM) stack, it is possible to uniquely reconstruct $X$ from the category $\mathrm{QCoh}(X)$. This is precisely the content of Tannaka reconstruction theorem, see \cite[Chapter 9]{SAG}. An analogous result holds in the abelian ($1$-categorical) setting, where it is known as the Rosenberg reconstruction theorem \cite{brandenburg2013rosenberg}.

In this paper, we focus on the study of two particularly well-behaved subcategories of $\mathrm{QCoh}(X)$, namely the category of perfect complexes and the category of bounded coherent sheaves. Perfect complexes and coherent sheaves can be defined over any ringed space in the following way:

\begin{definition}
 Let $(X,\m{O}_{X})$ be a ringed space. We say that an object $\m{F}$ in the derived category of $\m{O}_{X}$-modules is perfect if it is locally quasi-isomorphic to a bounded-below complex with coherent cohomology (i.e., pseudo-perfect) and finite Tor-amplitude. We say that $\m{F}$ is bounded and coherent if it is quasi-isomorphic to a bounded complex with coherent cohomology.	
\end{definition}

In particular, we study the $\infty$-categorical enhancements of these categories, which we denote by $\Perf(-)$ and $\Cohb(-)$, respectively. These enhancements have been explored in the context of algebraic and formal geometry in \cite{SAG}, as well as in non-Archimedean geometry in \cite{derhom,mathew}.

In algebraic geometry, the category of perfect complexes has been extensively studied, especially in relation to the K-theory of schemes, beginning with the seminal paper \cite{thomason2007higher}. This category possesses a fundamental property: in \cite{thomason2007higher}, it was shown that for a quasi-compact and quasi-separated scheme $X$, the category of perfect complexes, $\Perf(X)$, is the subcategory of compact generators of $\mathrm{QCoh}(X)$. Consequently, $\Perf(X)$ encodes rich geometric information about the scheme $X$.

Many geometric properties of an algebraic variety can be studied by examining its category of perfect complexes. A particularly useful tool is the ``non-commutative'' notion of smooth and proper categories, which has been extensively studied, originating with Kontsevich. In \cite{orlov,lunts2010categorical}, it was shown that an algebraic variety is smooth and proper if and only if its category of perfect complexes is also smooth and proper. Generalizing a theorem of Lunts \cite{lunts2010categorical}, Efimov proved in \cite{efimov} that the bounded derived category of coherent sheaves is of finite type (and thus smooth) even if the underlying algebraic variety is not smooth. Additionally, \cite{kontsevich1998triangulated} initiated the study of categories sharing properties with $\mathrm{Perf}(X)$, viewing these categories as ``non-commutative algebraic spaces.'' These and related ideas have been deployed with great success to prove theorems in algebraic geometry. For instance, van den Bergh, in \cite{van2004three}, uses non-commutative algebras to provide a proof of the Bondal-Orlov conjecture, showing that birational smooth projective 3-folds are derived equivalent.

The main aim of this paper is to study the category of perfect complexes on a rigid analytic variety and on a formal scheme, along with the notions of smooth and proper categories. We will demonstrate that in both settings, the situation differs significantly from the algebraic setting.
To develop these ideas in the formal and analytic settings, we encounter various difficulties. First, in derived analytic geometry, the definition of quasi-coherent sheaves is not straightforward. The recent work of Clausen and Scholze addresses this problem by developing the theory of \emph{Condensed Mathematics} \cite{clausen2019lectures} and using the \cat of solid quasi-coherent sheaves (see \cite{analytic,complex}). An alternative approach can be found in \cite{ben2024perspective,bambozzi2016dagger}, where the authors develop the theory of \emph{Bornological Mathematics} and propose defining the \cat of ``analytic quasi-coherent sheaves'' as the \cat of Ind-Banach modules. 

For our purposes, it suffices to work with the category of perfect complexes on an analytic variety. This provides a much simpler framework where the subtleties inherent in developing a full theory of ``analytic quasi-coherent sheaves'' are avoided. In particular, we define perfect complexes as in \cite[Definition 7.5]{derhom}. We utilize the foundations of derived non-Archimedean analytic geometry as developed in \cite{anmap,derhom,derived,porta2016higher}, and derived formal geometry as developed in Chapter 8 of \cite{SAG} and \cite{antonio2019derived}. We remark that the category of perfect complexes on an analytic variety can also be studied within the settings of Condensed and Bornological mathematics: indeed, as explained in \cite{andreychev2021pseudocoherent}, perfect complexes can be recovered as a full subcategory of the \cates of ``analytic quasi-coherent sheaves.'' However, unlike the usual setting of algebraic geometry, perfect complexes are no longer compact generators of these $\infty$-categories. We will show that the \cat of perfect complexes in the non-Archimedean analytic and formal settings is not as well-behaved as it is in the algebraic setting.

Throughout this paper, we consider rigid analytic varieties over a complete non-Archimedean field $k$ of characteristic $0$, and we assume all such varieties are locally of finite type. A standard reference for this topic is \cite{bosch2014lectures}.
Furthermore, \cite[Section 5.4]{bosch2014lectures} describes the analytification functor $(-)^{\mathrm{an}}$ from the category of schemes locally of finite type over $k$ to the category of rigid analytic varieties over $k$. In our main result, we characterize the essential image of this functor in certain favorable cases via the \cat of perfect complexes. Indeed, our main theorem is the following:

\begin{theorem}[\Cref{thm: main}]
    Let $X$ be a smooth and proper, quasi-compact, and separated rigid analytic variety over $k$. Then $\Perf(X)$ is smooth and proper as a $k$-linear \cat if and only if $X$ is algebraizable (i.e., there exists an algebraic space $Y$ such that $Y^{\mathrm{an}} \simeq X$).
\end{theorem}

From this theorem, we deduce that, in general, the category of perfect complexes of a smooth and proper rigid analytic variety is not smooth and proper, standing in stark contrast to the situation in algebraic geometry \cite{orlov}. This suggests that the notion of smooth and proper categories is intrinsically algebraic, and is not inherently well-adapted to the non-Archimedean and formal settings.

One of the main technical ingredients in our argument is an extension of the notion of simple objects from \cite{modg} to derived rings and rigid analytic varieties. In \cite{modg}, To\"en and Vaqui\'e define, for a compactly generated $k$-linear dg-category $\mathscr{C}$ of finite type, the derived moduli stack of objects 
\[
\mathscr{M}_{\mathscr{C}} \colon \mathrm{CAlg}_{k} \to \mathscr{S}
\]
from the \cat of derived commutative algebras to the \cat of spaces. The $\infty$-categorical analogue of this functor is informally described by sending a derived commutative algebra $A$ over $k$ to the maximal groupoid of the \cat of $k$-linear stable functors $\Fun^{\mathrm{st}}_{k}((\mathscr{C}^{\omega})^{\mathrm{op}}, \Perf(A))$.
In \cite{modg}, the authors also introduce the (underived) substack of simple objects $t_{0}\mathscr{M}_{\mathscr{C}}^{\mathrm{sim}} \subset t_{0}\mathscr{M}_{\mathscr{C}}$.

We extend both constructions to the derived analytic setting by introducing the (rigid) analytic derived stack of objects $\mathscr{AnM}_{\mathscr{C}}$, where $\mathscr{C}$ is a smooth and proper $\infty$-category, and the analytic derived substack of simple objects $\mathscr{AnM}_{\mathscr{C}}^{\mathrm{sim}}$. In \Cref{def: simple}, we also introduce a general notion of simple objects in a linear \cat, which allows us to construct the derived stack of simple objects $\m{M}^{\mathrm{sim}}_{\m{C}}$. Furthermore, within this general framework, we introduce the analytic derived stack of simple perfect complexes $\uline{\Perf}_{X}^{\mathrm{sim}}$ for a rigid analytic variety $X$.

We then establish an equivalence between the analytic derived stack of simple objects $\mathscr{AnM}_{\mathscr{C}}^{\mathrm{sim}}$ and the analytification of the algebraic stack of simple objects $(\mathscr{M}_{\mathscr{C}}^{\mathrm{sim}})^{\mathrm{an}}$ (\Cref{prop: MB sim vs AnMB sim}). Finally, if $X$ is a proper rigid analytic variety such that $\Perf(X)$ is smooth and proper, we show that $\uline{\Perf}^{\mathrm{sim}}_{X}$ is equivalent to the analytic derived stack $\mathscr{AnM}_{\Perf(X)}^{\mathrm{sim}}$ (\Cref{prop: Perf vs AnM}). These technical developments are detailed in Section 2 and constitute a key ingredient in the proof of our main theorem. Indeed, by employing the theory of analytic stacks developed in Section 2, we show in the final section that there exists an embedding of $X$ into the coarse moduli space of $(t_0\mathscr{M}_{\Perf(X)}^{\mathrm{sim}})^{\mathrm{an}}$ (see \Cref{prop: X open into the analytification}). This embedding, together with the concept of non-Archimedean Moishezon spaces introduced in \cite{conrad}, provides the main tools used in the proof of \Cref{thm: main}.

As a consequence of our main theorem, we prove two related results on algebraizability for formal schemes. In this setting, we consider formal schemes over $\Spf(\m{O}_{k})$, where $\m{O}_{k}$ is a DVR with maximal ideal $\mathfrak{m}$ and fraction field $k$. We assume that $\m{O}_k$ is complete with respect to the $\mathfrak{m}$-adic topology. All formal schemes are assumed to be separated, quasi-compact, topologically almost of finite presentation, and locally Noetherian.

\begin{theorem}[\Cref{thm: main formal setting}]
    Let $\F{X}$ be a smooth and proper formal scheme over $\Spf(\m{O}_{k})$. Then $\Perf(\F{X})$ is smooth and proper as an $\m{O}_k$-linear stable \cat if and only if $\F{X}$ is algebraizable. 
\end{theorem} 

\begin{theorem}[\Cref{thm: algebraization via coherent categories}]
    Let $\F{X}$ be a proper formal scheme such that $\Cohb(\F{X})$ is smooth as a \cat in $\catperf_{\m{O}_{k}}$, and $\F{X}^{\rig}$ is smooth. Then $\F{X}$ is algebraizable. 
\end{theorem}

These theorems are, of course, closely related. We prove both theorems by using \Cref{prop: algebrization} to reduce these statements to \Cref{thm: main}.

\Cref{thm: algebraization via coherent categories} shows that, in general, for a proper non-algebraizable formal scheme $\F{X}$, we cannot expect $\Cohb(\F{X})$ to be smooth. On the other hand, for a separated scheme $X$ over a perfect field, Lunts demonstrated in \cite{lunts2010categorical} that $\Cohb(X)$ is smooth. More generally, Efimov showed in \cite{efimov} that under mild conditions on $X$, the category $\Cohb(X)$ is of finite type. Since being of finite type implies smoothness, the above theorem confirms that these classical algebraic statements cannot be broadly generalized to the setting of formal schemes.
\subsubsection*{Relation to other works.} 
\Cref{thm: main} is a generalization of the main result of \cite{tova} to the non-Archimedean setting. In \cite{tova}, the authors prove  that a complex analytic, connected, smooth, and proper variety is algebraizable if and only if its category of perfect complexes is smooth and proper as a $\mathbb{C}$-linear DG-category. We remark that our proof strategy is modelled after the arguments in  \cite{tova}, although we have to deal with many new subtleties that are absent in the complex setting. 
We can interpret our result and the results from\cite{tova} as evidence of the fact that the category of perfect complexes in the analytic and formal settings is not as well-behaved as in the algebraic setting. The key issue is that the definition of smooth and proper categories is not well adapted to the analytic setting. In future work we intend to return to the problem of developing a notion of smooth and proper category which is better suited to analytic geometry.

\subsubsection*{Acknowledgments} I wish to express my gratitude to my advisors Mauro Porta and Nicol\`o Sibilla for the numerous conversations we've had. Many of the ideas presented in this paper have stemmed from these discussions; I also extend my thanks to Emanuele Pavia for the valuable conversations for many suggestions that have improved the paper's redability. I would also like to
thank Andrea Ricolfi and Enrico Lampetti for valuable suggestions on an earlier version of this paper.

\subsubsection*{Notations and conventions}
In this paper, we work with a fixed base field $k$ of characteristic $0$. In the context of formal schemes, we further assume that $k$ is the field of fractions of a discrete valuation ring (DVR) $\mathscr{O}_{k}$, complete with respect to the $\mathfrak{m}$-adic topology. Additionally, in the non-Archimedean setting, we require the field $k$ to be complete. Throughout this paper, we freely employ the language of derived geometry, $\infty$-categories, and higher algebra. In this setting, we follow the conventions of \cite{HTT}; in particular, we denote by $\m{S}$ the $\infty$-category of spaces. All rings and modules considered are derived unless explicitly stated otherwise. In particular, a module over an algebra $A$ will be an object in its derived $\infty$-category, as explained in \cite[Chapter 7]{HA}. We denote this category by $\Mod_{A}$, its subcategory of algebra objects by $\CAlg_{A}$, and the opposite category of $\CAlg_{A}$ by $\Aff_{A}$. When discussing (derived) algebraic stacks, we equip the category $\Aff_{k}$ (or $\Daff_{k}$) with the \'etale topology. In the rigid analytic setting, we denote by $\Dafd_{k}$ the \cat of derived affinoid spaces, as defined in \cite[Definition 7.3]{derived}. Similarly, in the underived setting, we denote by $\Afd_{k}$ the category of affinoid spaces over $k$. When working with derived analytic stacks, we consider the \'etale topology on the category $\Dafd_{k}$. For derived geometry, we follow the notation of \cite{SAG} in the algebraic and formal settings, while in the rigid analytic setting, our main references are \cite{anmap} and \cite{derhom}. Furthermore, all formal schemes and rigid analytic varieties are assumed to be topologically of finite type, quasi-compact, and separated.

\section{Moduli of objects for \texorpdfstring{$k$}{k}-linear stable \texorpdfstring{\cates}{oo-categories}}
In this section, we revisit the construction of the moduli stack of objects in a DG category as described in \cite{modg}, adopting the language of $\infty$-categories. Since we shall employ this construction both for \cates in $\Pr^{\mathrm{L},\omega}$ and for small \cates, we begin with some remarks on stable Morita theory, which facilitates the transition between these two settings.

\subsection{A perspective on Morita theory via Ind-completion}\label{sec1}


Here we describe Morita theory for stable $\infty$-categories. Building upon unstable Morita theory, we show how to extend this correspondence to the stable setting. We denote by $\mathrm{Cat}^{\Perf,\kappa}$, where $\kappa$ is a regular cardinal, the \cat that has as objects all small \cates with $\kappa$-small colimits, which are idempotent complete and whose morphisms preserve $\kappa$-small filtered colimits (see \cite{HTT} and \cite{blumberg2013universal}). Additionally, $\Pr^{\mathrm{L},\kappa}$ denotes the \cat of presentable \cates generated by $\kappa$-compact objects, where morphisms preserve $\kappa$-compact objects. 

The functor $\Ind_{\kappa}(-)$ is the formal completion with respect to $\kappa$-filtered colimits; the functor $(-)^{\kappa}$ sends a category to the full subcategory of its $\kappa$-compact objects. By the definition of an accessible \cat, this subcategory is small and $\kappa$-small cocomplete, as $\kappa$-small limits commute with $\kappa$-filtered colimits. Furthermore, it is idempotent complete, since compact objects are closed under retracts.

The main result in the unstable setting is the following theorem.

\begin{theorem}
    Let $\kappa$ be a regular cardinal. Then the following hold:
    \begin{enumerate}
        \item There is an equivalence of \cates given by the following pair of adjoint functors:
          \begin{equation*}
            \begin{tikzcd}
             \mathrm{Cat}^{\Perf,\kappa} \arrow[r, shift left=1ex, "\mathrm{Ind}_{\kappa}"{name=G}] & \Pr^{\mathrm{L},\kappa}\arrow[l, shift left=.5ex, "(-)^{\kappa}"{name=F}]
            \arrow[phantom, from=F, to=G, , "\scriptscriptstyle\boldsymbol{\top}"].
            \end{tikzcd}
        \end{equation*}
        \item The \cat $\mathrm{Pr}^{\mathrm{L},\kappa}$ is presentable.
        \item The obvious functor $\mathrm{Pr}^{\mathrm{L},\kappa} \to \pr $ commutes with small colimits.
    \end{enumerate}
\end{theorem}

\begin{proof}
    The first statement is exactly \cite[Proposition 5.5.7.8]{HTT}, and the second follows from \cite[Proposition 2.3]{aoki2025higherpresentablecategorieslimits}. The third statement can be deduced using \cite[Proposition 5.5.7.6, and Theorem 5.5.3.18]{HTT} and by considering that, from \cite[Corollary 5.5.3.4]{HTT}, we have an equivalence 
    \[
    (\Pr\nolimits^{\mathrm{R}})^{\mathrm{op}} \simeq \Pr\nolimits^{\mathrm{L}}.
    \]
    Using \cite[Proposition 5.5.7.2]{HTT}, we can deduce that this equivalence restricts to an equivalence
    \[
    (\Pr\nolimits^{\mathrm{R,\kappa}})^{\mathrm{op}} \simeq \Pr\nolimits^{\mathrm{L}, \kappa}.
    \]
\end{proof}


We can apply the first part of the theorem to the case where $\kappa=\omega$. In this setting, we denote the \cat $\mathrm{Cat}^{\Perf,\omega}$ simply by $\mathrm{Cat}^{\Perf}$, which is the \cat of small, idempotent complete \cates with finite colimits. This gives us the following equivalence:
\begin{equation*}
    \catperf \simeq \promega.
\end{equation*}

By restricting to the subcategory of stable objects $\catperf_{\mathrm{st}}$ with exact functors, we obtain the following version of stable Morita theory:

\begin{proposition}\label{prop:stable_Morita} 
    There is an equivalence of $\infty$-categories
    \begin{equation*}
    \catperf_{\mathrm{st}} \simeq \promegast.
    \end{equation*}
\end{proposition}

\begin{proof}
    This essentially follows from the above equivalence, noting that stable \cates are closed under $\Ind_{\omega}$ (usually denoted simply by $\Ind$) \cite[Proposition 1.1.3.6]{HA}. Additionally, the subcategory of compact generators of a stable \cat is itself stable.
\end{proof}

We note that we can also start with a stable, small, \cat $\mathscr{C}$ that is not idempotent complete and obtain the same result as above by replacing $\mathscr{C}$ with $\mathscr{C}^{\mathrm{idem}}$, here, as explained in \cite{HTT}, $\mathscr{C}^{\mathrm{idem}}$ denotes the idempotent completion of $\mathscr{C}$. Moreover, $\mathscr{C}^{\mathrm{idem}}$ remains small \cite[Proposition 5.4.2.4]{HTT} and stable \cite[Corollary 1.1.3.7]{HA}.
We will now extend the above stable Morita theory to $k$-linear $\infty$-categories. We begin by noting that, as described in \cite[Theorem 3.1]{blumberg2013universal}, there is a symmetric monoidal structure on $\mathrm{Cat}^{\Perf}_{\mathrm{st}}$. This provides the context for the following definition.

\begin{definition}
    We define the \cat of stable, idempotent complete, $k$-linear \cates as 
    \begin{equation*}
        \catperf_{k} \coloneqq \f{Mod}_{\f{Perf}(k)}{\catperf_{\mathrm{st}}} .
    \end{equation*}
    The \cat of stable, presentable, compactly generated $k$-linear \cates is:
    \begin{equation*}
        \promega_{k} \coloneqq \f{Mod}_{\f{Mod}_{k}}{\promegast} .
    \end{equation*}
\end{definition}

\begin{remark}
Note that it is natural to consider, in $\mathrm{Cat}^{\Perf}$, modules over $\mathrm{Perf}(k)$ rather than over $\mathrm{Mod}_{k}$, since $\mathrm{Mod}_{k}$ is not a small category, and in $\mathrm{Cat}_{\infty,k}$ we only deal with small $\infty$-categories. However, under the equivalence from \Cref{prop:stable_Morita}, $\mathrm{Mod}_{k}$ corresponds precisely to $\mathrm{Perf}(k)$.
\end{remark}

\begin{lemma}\label{lem:symmetric_monoidal_functors_and_modules}
Let $\mathscr{C}$ and $\mathscr{D}$ be two symmetric monoidal \cates, and let $F \colon \mathscr{C} \to \mathscr{D}$ be a symmetric monoidal equivalence. If $A$ is an algebra object in $\mathscr{C}$, then $F$ induces an equivalence:
\begin{equation*}
    \f{Mod}_{A}\mathscr{C} \simeq \f{Mod}_{F(A)}\mathscr{D} .
\end{equation*}
\end{lemma}
\begin{proof}
The proof is straightforward.
\end{proof}

From the above lemma we can deduce the following characterization for the $k$-linear stable Morita \cat.
\begin{proposition}
The pair of adjoint functors
          \begin{equation*}
            \begin{tikzcd}
            \mathrm{Cat}^{\Perf}_{k} \arrow[r, shift left=1ex, "\Ind_{\omega}"{name=G}] & \Pr^{\mathrm{L},\omega}_{k}\arrow[l, shift left=.5ex, "(-)^{\omega}"{name=F}]
            \arrow[phantom, from=F, to=G, , "\scriptscriptstyle\boldsymbol{\top}"]
            \end{tikzcd}
        \end{equation*}
        forms an equivalence.
\end{proposition}

\subsection{Smooth and proper \cates and algebras}
In this subsection, we begin by reviewing several established results concerning smooth and proper algebras and $\infty$-categories.

Let $A$ be an $\mathbb{E}_{\infty}$-ring and $B$ an $\mathbb{E}_{1}$-algebra in $\Mod_{A}$. 
\begin{definition}
    We say that $B$ is \textit{proper} over $k$ if it is a perfect complex in $\f{Mod}_{A}$.
\end{definition}

\begin{definition}
    We say that $B$ is \textit{smooth} over $k$ if it is a compact object in $\f{Mod}_{B \otimes_{k} B^{\mathrm{op}}}$.
\end{definition}

We have the following useful characterization:
\begin{proposition}[{\cite[4.6.4.4 and 4.6.4.12]{HA}}]\label{prop:caracterization of smooth and proper}\leavevmode 

\begin{enumerate}
    \item An algebra $B$ is proper if and only if for every $M \in \LMod_{B}$, if $M$ is dualizable in $\LMod_{B}(\Mod_{A})$, then $M$ is dualizable in $\Mod_{A}$.
    \item An algebra $B$ is smooth if and only if for every $M \in \LMod_{B}(\Mod_{A})$, if $M$ is dualizable as an $A$-module, then $M$ is dualizable in $\LMod_{B}(\Mod_{A})$.  
\end{enumerate}
\end{proposition}
\begin{lemma}\label{lem: perfects on a smooth and proper algebra}
    If $B$ is smooth and proper, then $\LMod_{B}(\Perf(A)) \simeq \Perf(B)$.
\end{lemma}
\begin{proof}
    This follows directly from the characterization of smooth and proper algebras, noting that dualizable objects coincide with perfect complexes in $\Mod_{A}$.
\end{proof}

In the setting of presentable, compactly generated \cates, we define smoothness and properness as follows. Let $\m{C}$ be an \cat in $\mathrm{Pr^{L}_{k}}$, and let $\m{C}^{\vee} = \Fun^{\mathrm{L}}_{k}(\m{C}, \Mod_{k})$ be its dual category, equipped with the evaluation map $\mathrm{ev} \colon \m{C} \otimes \m{C}^{\vee} \to \Mod_{k}$.

\begin{definition}
    A \cat $\m{C} \in \mathrm{Pr}^{L}_{A}$ is \textit{dualizable} if there exists a coevaluation map
    \begin{equation*}
        \mathrm{coev} \colon \Mod_{A} \to \m{C} \otimes \m{C}^{\vee}
    \end{equation*}
    classifying $\m{C}$ as a module over $\m{C} \otimes \m{C}^{\vee}$, such that both compositions induced by evaluation and coevaluation: 
    \begin{equation*}
        \begin{tikzcd}
            {\m{C}} && {\m{C} \otimes \m{C}^{\vee} \otimes \m{C}} && {\m{C}}
            \arrow["{\mathrm{id} \otimes \mathrm{coev}}", from=1-1, to=1-3]
            \arrow["{\mathrm{ev} \otimes \mathrm{id}}", from=1-3, to=1-5]
        \end{tikzcd}
    \end{equation*}
    \begin{equation*}
        \begin{tikzcd}
            {\m{C}^{\vee}} && {\m{C}^{\vee} \otimes \m{C} \otimes \m{C}^{\vee}} && {\m{C}^{\vee}}
            \arrow["{\mathrm{coev} \otimes \mathrm{id}}", from=1-1, to=1-3]
            \arrow["{\mathrm{id} \otimes \mathrm{ev}}", from=1-3, to=1-5]
        \end{tikzcd}
    \end{equation*}
    are homotopic to the identity.
\end{definition}

\begin{definition}\label{def: proper cat}
    We say that $\m{C} \in \mathrm{Pr}^{\mathrm{L},\omega}_{A}$ is \textit{proper} if the evaluation map is a morphism in $\mathrm{Pr}^{\mathrm{L},\omega}_{A}$ (i.e., it admits a right adjoint and preserves compact objects).
\end{definition}

\begin{definition}\label{def: smooth cat}
    We say that $\m{C} \in \mathrm{Pr}^{\mathrm{L},\omega}_{A}$ is \textit{smooth} if it is dualizable and the coevaluation map is a morphism in $\mathrm{Pr}^{L,\omega}_{A}$.
\end{definition}

\begin{definition}
    We say that $\m{C} \in \mathrm{Pr}^{\mathrm{L},\omega}_{A}$ is of \textit{finite type} if there exists a compact $A$-algebra $C$ such that $\m{C}$ is isomorphic to $\Mod_{C}$.
\end{definition}

We also define the notion of smoothness (respectively, properness) for a small, $A$-linear, idempotent complete \cat by requiring that its idempotent completion is smooth (respectively, proper) as in \Cref{def: smooth cat} (respectively, \Cref{def: proper cat}). In this context, a small, $k$-linear, idempotent complete \cat $\cB$ is proper if and only if, for every pair of objects $b, b'$ in $\cB$, the $A$-module $\mathrm{Hom}_{\cB}(b,b')$ is a perfect complex.

\begin{lemma}\label{lem: compact generator}
    Let $\m{C}$ be a smooth $A$-linear \cat in $\mathrm{Pr}^{\mathrm{L},\omega}_{A}$; then $\m{C}$ admits a compact generator.
\end{lemma}
\begin{proof}
    See the proof of \cite[Proposition 11.3.2.4]{SAG}.
\end{proof}

\begin{corollary}\label{cor: categories vs algebras}
    An $A$-linear \cat in $\mathrm{Pr}^{\mathrm{L},\omega}_{A}$ is smooth and proper if and only if it is equivalent to $\LMod_{B}$ for some smooth and proper algebra $B$. In this setting, $B$ can be taken to be the algebra of endomorphisms of a compact generator.
\end{corollary}
\begin{proof}
    If $\mathscr{C}$ is a smooth and proper $A$-linear presentable \cat, using \cite[Theorem 7.1.2.1]{HA} and \Cref{lem: compact generator}, we can deduce that there exists an algebra $B$ such that $\mathscr{C}$ is equivalent to $\LMod_{B}$. The statement is then clear, as the categorical definitions of smoothness and properness directly correspond to their algebraic counterparts.
\end{proof}

\begin{lemma}[{\cite[Lemma 2.6]{modg}}]\label{lem: smooth and proper small categories}
    An \cat $\m{C}$ in $\mathrm{Pr}^{\mathrm{L},\omega}_{A}$ is smooth and proper if and only if its subcategory of compact objects $\m{C}^{\omega}$ is smooth and proper. 
\end{lemma}

\subsection{The stack of pseudo-perfect objects}
Here, we review the main construction of \cite{modg} using the language of $\infty$-categories. Although the DG category setting is used in \cite{modg}, similar ideas can be applied in the $\infty$-categorical setting, as demonstrated in \cite{antieau2014brauer}. For a comprehensive overview of this construction in the language of $\infty$-categoires, see also \cite{lampetti2025goodmodulispaceconstructible}.

Let $k$ be a discrete field and let $\mathscr{C}$ be a stable, presentable, $k$-linear, compactly generated \cat (i.e., $\mathscr{C} \in \mathrm{Pr}^{\mathrm{L},\omega}_{k}$).

\begin{definition}
    A functor from $\mathscr{C}$ to $\Mod_{A}$ is called \emph{pseudo-perfect} if it sends compact objects of $\mathscr{C}$ to perfect $A$-modules.
\end{definition}

The following functor, $\mathscr{M}_{\mathscr{C}}$, was constructed in \cite{modg} and can be described as the $\infty$-functor: 
\begin{equation*}
\begin{tikzcd}[row sep=0.01in]
\mathscr{M}_{\mathscr{C}} \colon \mathrm{CAlg}_{k} \arrow[r] & \mathscr{S}\\
A \arrow[r, maps to]                               & {(\mathrm{Fun}^{\mathrm{R}}_{\mathrm{ps},k}(\mathscr{C}^{\mathrm{op}},\mathrm{\Mod}_{A}))^{\simeq}}
\end{tikzcd}
\end{equation*}
where the subscript ``ps'' denotes pseudo-perfect functors.

Informally, the functor
\begin{equation*}
    \Mod_{-}\colon \mathrm{CAlg}_{k} \to \pr
\end{equation*}
is defined by sending a morphism $A\to B$ to the base change functor $-\otimes_{A}B \colon \Mod_{A} \to \Mod_{B}$. We then use the functoriality of $\Fun^{\mathrm{R}}_{\mathrm{ps},k}(\mathscr{C}^{\mathrm{op}},-)$, which relies on the fact that $\Fun^{\mathrm{R}}_{k}(-,-)$ is already known to be a functor, and that the tensor product preserves perfect objects.

We can also describe $\Fun^{\mathrm{R}}_{\mathrm{ps},k}(\mathscr{C}^{\mathrm{op}},\Mod_{A})$ equivalently as follows:
\begin{equation}
\begin{aligned}
    \mathrm{Fun}^{\mathrm{R}}_{\mathrm{ps},k}(\mathscr{C}^{\mathrm{op}},\f{Mod}_{A}) & \simeq \f{Fun}^{\mathrm{R},\omega}_{k}(\mathscr{C}^{\mathrm{op}},\f{Mod}_{A}) \\
    & \simeq \f{Fun}_{k}^{\mathrm{R},\omega}((\mathrm{Ind}(\mathscr{C}^{\omega}))^{\mathrm{op}},\f{Mod}_{A}) \\ 
    & \simeq \f{Fun}^{\mathrm{L},\omega}_{k}(\mathrm{Ind}(\mathscr{C}^{\omega}),\f{Mod}_{A}^{\mathrm{op}})^{\mathrm{op}}\footnotemark \\
    & \simeq \f{Fun}_{k}^{\mathrm{st}}(\mathscr{C}^{\omega}, \f{Perf}(A)^{\mathrm{op}})^{\mathrm{op}} \\
    & \simeq \Fun^{\mathrm{st}}_{k}((\mathscr{C}^{\omega})^{\mathrm{op}}, \Perf(A)).
\end{aligned}
\end{equation}
The last equivalence provides the $\infty$-categorical analog of the one used in \cite{modg}.

\footnotetext{This follows from the universal property of Ind-completion (i.e., the \cat of functors preserving filtered colimits from $\mathrm{Ind}(\mathscr{C})$ to $\mathscr{D}$ is equivalent to the \cat of functors from $\mathscr{C}$ to $\mathscr{D}$) and the fact that all colimits can be obtained as a composition of finite colimits and filtered colimits.}

\begin{remark}\label{rmk: MB vs LModB} 
We observe that the \cat $\mathscr{M}_{\LMod_{B}}(A)$ is equivalent to $\LMod_{B}(\Perf(A))$. Indeed, a pseudo-perfect, $k$-linear functor from $\LMod_{B}^{\mathrm{op}}$ to $\Mod_{A}$ corresponds to a perfect complex $F$ on $A$ and a map 
\[
\mathrm{End}_{\LMod_{B}^{\mathrm{op}}}(B,B) \to \mathrm{End}_{\Mod_{A}}(F).
\]
It suffices to observe that, as an $\mathbb{E}_{1}$-algebra, $\mathrm{End}_{\LMod_{B}^{\mathrm{op}}}(B,B)$ is equivalent to $B$.
\end{remark}

\begin{theorem}[{\cite{modg}}]
    The functor $\mathscr{M}_{\mathscr{C}} \colon \mathrm{CAlg}_{k} \to \mathscr{S}$ is a derived stack, and if $\mathscr{C}$ is of finite type, $\mathscr{M}_{\mathscr{C}}$ is locally geometric.
\end{theorem}

\begin{definition}
    We call the functor $\mathscr{M}_{\mathscr{C}}$ the moduli stack of objects of $\mathscr{C}$. 
\end{definition}

Now we can view $\mathscr{M}_{-}$ as a functor from $(\promega_{k})^{\mathrm{op}} \to \mathrm{dSt}$ (where $\mathrm{dSt}$ denotes the \cat of derived stacks).

\section{Simple objects in an \cat and moduli of objects in rigid geometry}

In this section, we define several analytic stacks. In particular, we construct a rigid analytic version of the derived moduli stack of objects introduced in \cite{modg}. Given a smooth and proper $k$-linear \cat $\cC$, we define the analytic stack of objects in $\cC$, denoted by $\mathscr{AnM}_{\cC}$. We adapt the notion of simple objects from \cite{modg} to define an analytic variant of the moduli stack of simple objects, as in \cite[Definition 3.21]{modg}. Specifically, for a smooth and proper $k$-linear \cat $\cC$, we construct the analytic substack $\mathscr{AnM}_{\cC}^{\text{sim}} \subset \mathscr{AnM}_{\cC}$. Furthermore, we introduce the analytic stack of simple perfect complexes over a rigid analytic variety $X$, denoted by $\uline{\Perf}^{\text{sim}}_{X}$. The goal of this section is to compare these analytic stacks with the analytification of the algebraic stack $\m{M}^{\mathrm{sim}}_{\Perf(X)}$. In particular, we will prove the following result, which plays a crucial role in the proof of \Cref{thm: main}.

\begin{theorem}
    Let $X$ be a proper, quasi-compact, and separated rigid analytic variety over a non-Archimedean field $k$, such that $\Perf(X)$ is a smooth and proper $k$-linear \cat. Then, there is an equivalence of analytic stacks between $\uline{\Perf}_{X}$ and the analytification of the algebraic stack $\mathscr{M}_{\Perf(X)}$. Moreover, this equivalence restricts to simple objects, and we obtain an equivalence of $1$-stacks between $t_0\uline{\Perf}^{\mathrm{sim}}_{X}$ and $(t_0\m{M}_{\Perf(X)}^{\mathrm{sim}})^{\mathrm{an}}$.
\end{theorem}

We obtain this theorem as a consequence of \Cref{prop: MB sim vs AnMB sim}, \Cref{prop: Perf vs AnM}, and \Cref{cor: analytifiction of the trouncation}.

We begin by defining the general notion of a simple object in \Cref{Sec2.1}. Then, in \Cref{sec2.2}, \Cref{sec2.3}, and \Cref{sec2.4}, we introduce the main stacks considered in this paper and adapt the concept of simple objects to this setting. Finally, we prove \Cref{prop: MB sim vs AnMB sim} in \Cref{secComparison1}, followed by the proofs of \Cref{prop: Perf vs AnM} and \Cref{cor: analytifiction of the trouncation} in \Cref{secComparison2}.

We denote by $\Afd_{k}$ the category of affinoid spaces (as defined in \cite{bosch2014lectures}) over a fixed complete non-Archimedean field $k$, of characteristic $0$, endowed with a non-trivial valuation. Affinoid spaces play the same role in rigid geometry as affine schemes do in algebraic geometry. In our context, the most relevant examples of affinoid spaces for constructing rigid analytic varieties are those whose algebra of functions is a quotient of the Tate algebra $k\langle t_{1}, \dots, t_{n} \rangle$ by a finitely generated ideal; we will refer to these as affinoid spaces of finite type. For precise definitions of affinoid spaces and rigid analytic varieties, we refer to \cite{bosch2014lectures}. We implicitly assume that all rigid analytic varieties are quasi-compact, separated, and of finite type.

To properly formulate the notion of derived stacks in rigid analytic geometry, we adopt the framework developed in \cite{derhom, derived, porta2016higher}. Throughout our discussion of stacks, we consistently work with the \'etale topology.

We denote by $\Dafd_{k}$ the \cat of derived affinoid spaces over $k$, as described in \cite[Definition 7.3]{derived}. Additionally, we denote by $\Sp(A)$ the underived affinoid space associated with an underived affinoid algebra $A$. It is important to note that, following the conventions of \cite{derhom}, derived affinoids do not strictly take the form $\Sp(A)$ as they do in the classical setting. However, if $S$ is a derived affinoid, its algebra of functions is a derived affinoid algebra (i.e., a connective $\mathbb{E}_{\infty}$-algebra $A$ in $\Mod_{k}$ such that $\pi_{0}(A)$ is an underived affinoid algebra and the $\pi_{i}(A)$ are finitely generated modules over $\pi_{0}(A)$).

\subsection{Simple objects in \cates}\label{Sec2.1}
Let $\m{V}$ be a commutative algebra object in $\prst$ equipped with a t-structure, and let $\m{C}$ be a presentable $\cV$-linear \cat (i.e., an object in $\Mod_{\m{V}}(\pr)$). Since $\m{C}$ is a module over $\m{V}$, we denote by
\begin{equation*}
    \m{V}\otimes \m{C} \to \m{C}
\end{equation*}
the canonical action of $\m{V}$ on $\m{C}$. Let $F$ be an object in $\m{C}$. Since the functor $-\otimes F$ commutes with colimits, it admits a right adjoint, which we denote by $\f{Hom}_{\m{V}\otimes \m{C}}(F,-)$.
\begin{equation}\label{eq: external hom}
\begin{tikzcd}
	{\mathscr{C}} &&& {\mathscr{V}}
	\arrow[""{name=0, anchor=center, inner sep=0}, "{\Hom_{\mathscr{V}\otimes \mathscr{C}}(F,-)}"', shift right=2, from=1-1, to=1-4]
	\arrow[""{name=1, anchor=center, inner sep=0}, "{-\otimes F}"', shift right=2, from=1-4, to=1-1]
	\arrow["\dashv"{anchor=center, rotate=-90}, draw=none, from=1, to=0]
\end{tikzcd}
\end{equation}

\begin{remark}\label{rmk:map simple}
Note that for every object $F \in \m{C}$, the action of $\m{V}$ induces a map:
\begin{equation*}
    \begin{tikzcd}
\mathbf{1}_{\m{V}} \otimes F \arrow[r, "id"] & F
\end{tikzcd}
\end{equation*}
By adjunction, there is also a canonical map:
\begin{equation*}
    \mathbf{1}_{\m{V}} \to \f{Hom}_{\m{V}\otimes \m{C}}(F,F).
\end{equation*}
We now use this map to define simple objects, and in the subsequent subsections, we will apply this definition within various settings. Note that this map depends only on the action and the adjunction, rather than on the relative tensor product.
\end{remark}

\begin{remark}
Even in the noncommutative setting, within the \cat $\f{LMod}_{B}(\m{C})$ where $B$ is an associative algebra in a symmetric monoidal \cat $\m{C}$, it is possible to give a meaningful definition of dualizable object. The definitions and main properties are essentially the same as in the commutative setting, albeit with a few differences. For instance, if $F \in \f{LMod}_{B}$, then its left dual $F^{\vee}$, if it exists, belongs to $\f{RMod}_{B}$. Moreover, the relative tensor product $F^{\vee} \otimes_{B} F$ is an object of $\m{C}$ equipped with a canonical map from the monoidal unit $\mathbf{1}_{\m{C}}$, induced by the identity morphism. Indeed, by \cite[Proposition 4.6.2.18]{HA}, we have
\begin{equation*}
    \f{Hom}_{\mathrm{RMod}_B}(F^{\vee}\otimes \mathbf{1}_{\m{C}} , F^{\vee})\simeq \f{Hom}_{\m{C}}(\mathbf{1}_{\m{C}},F^{\vee}\otimes_{\m{C}} F).
\end{equation*}
Note that the first tensor product is associated with the action of $\mathscr{C}$ on $F$, whereas the latter is the exterior tensor product (i.e., $-\otimes_{B}-: \mathrm{RMod}_{B}({\mathscr{C}}) \times \mathrm{LMod}_{B}({\mathscr{C}}) \to \mathscr{C}$). For precise definitions and the main properties of the relative tensor product, we refer the reader to \cite[Section 4.6.2]{HA}.
\end{remark} 

\begin{definition}
    Let $\m{V}$ be a symmetric monoidal \cat equipped with a t-structure. We say that an object $C \in \m{V}$ has Tor-amplitude strictly smaller than zero if the homological amplitude of the tensor product $C\otimes D$ is strictly smaller than zero for every $D$ in the heart of $\m{V}$.
\end{definition}

\begin{definition}\label{def: simple}
    Let $\m{V}$ be a presentable symmetric monoidal \cat equipped with a t-structure, and let $\m{C}$ be a $\cV$-linear \cat. We say that an object $F$ in $\m{C}$ is \emph{simple} if 
    \begin{equation*}
        \cofib(\mathbf{1}_{\m{V}} \to \f{Hom}_{\m{V}\otimes \m{C}}(F,F))
    \end{equation*}
    has Tor-amplitude strictly smaller than zero.
\end{definition}

In the following subsections, we analyze this definition within the algebraic and rigid analytic settings. In particular, we will show that this general definition recovers the non-Archimedean analogues of the definitions of simple objects given in \cite{tova} for the stacks $\m{M}_{\mathscr{C}}$, $An\m{M}_{\mathscr{C}}$, and $\uline{\Perf}_{X}$, where $\mathscr{C}$ is a smooth and proper \cat and $X$ is a smooth and proper rigid analytic variety. Furthermore, using \Cref{cor: categories vs algebras}, we can reduce to the case where $\cC$ is equivalent to $\LMod_{B}$ for a smooth and proper algebra $B$. This is precisely the case we will begin studying.

\subsection{Simple objects in \texorpdfstring{$\m{M}_{B}$}{MB} with \texorpdfstring{$B$}{B} a smooth and proper \texorpdfstring{$\mathbb{E}_{1}$}{E1}-algebra}\label{sec2.2}\leavevmode \\

To define the substack of simple objects of $\m{M}_{B}$, we fix a connective $\mathbb{E}_{\infty}$-algebra $A$ in $\Mod_{k}$ and a smooth and proper $\mathbb{E}_{1}$-algebra $B$ in $\Mod_{k}$. Based on the general discussion above, we now define simple objects in $\m{M}_{B}(A)$. 

\begin{remark}\label{rmk: adjunction Mb}
We can apply \Cref{def: simple} by taking $\m{A}=\mathrm{Pr}^{\mathrm{L}}$, $\m{V}=\Mod_{A}$, and $\m{C}=\LMod_{B} \otimes \Mod_{A}$. Here, the action is informally given by sending $(M, M' \otimes N)$ to $M \otimes_{A} M' \otimes N$. This induces an adjunction:
\begin{equation*}
\begin{tikzcd}
	{\LMod_{B} \otimes \Mod_{A}} &&& {\Mod_{A}}
	\arrow[""{name=0, anchor=center, inner sep=0}, "{\Hom_{B \otimes A}(F,-)}"', shift right=2, from=1-1, to=1-4]
	\arrow[""{name=1, anchor=center, inner sep=0}, "{F \otimes_{A}-}"', shift right=2, from=1-4, to=1-1]
	\arrow["\dashv"{anchor=center, rotate=-90}, draw=none, from=1, to=0]
\end{tikzcd}
\end{equation*}
with $F \in \LMod_{B} \otimes \Mod_{A}$, where we denote by $\Hom_{B \otimes A}$ the external $\Hom$ constructed in \eqref{eq: external hom}.
\end{remark}

We will now prove that this induces an adjunction:
\begin{equation*}
    \begin{tikzcd}
	{\Perf(B \otimes A)} &&& {\Perf(A)}
	\arrow[""{name=0, anchor=center, inner sep=0}, "{\Hom_{B \otimes A}(F,-)}"', shift right=2, from=1-1, to=1-4]
	\arrow[""{name=1, anchor=center, inner sep=0}, "{F \otimes_{A}-}"', shift right=2, from=1-4, to=1-1]
	\arrow["\dashv"{anchor=center, rotate=-90}, draw=none, from=1, to=0]
\end{tikzcd}
\end{equation*}
with $F \in \Perf(B \otimes A)$. Note that the tensor product used here is the usual derived tensor product of $A$-modules, since the tensor product of perfect complexes is again perfect.

\begin{proposition}\label{prop:afjunction for simple objects algebra}
    Under the hypotheses described above, the functor
    \[
    F\otimes_{A}-\colon \Perf(A) \to \Perf(B\otimes A)
    \]
    admits a right adjoint, which we denote by $\Hom_{B \otimes A}(F,-)$.
\end{proposition}

\begin{proof}
    Consider the following diagram, in which the vertical maps are fully faithful. Note that at the level of module categories, we already know that the adjunction exists:
    \begin{equation*}
        \begin{tikzcd}
\LMod_{B} \otimes \Mod_{A} \arrow[rr, "{\Hom_{B \otimes A}(F,-)}"', shift right] &  & \Mod_{A} \arrow[ll, "F \otimes_{A} -"', shift right] \\
                                                                                &  &                                                 \\
\Perf(B \otimes A) \arrow[uu]                                        &  & \Perf(A) \arrow[uu] \arrow[ll, "F \otimes_{A}-"]    
\end{tikzcd}
    \end{equation*}
    Since the tensor product functors coincide, it suffices to prove that if $F, G \in \Perf(B \otimes A)$, then $\Hom_{B\otimes A}(F,G)$ remains a perfect complex.
To this end, let $\m{W}$ be the smallest stable full subcategory of $\LMod_{B}\otimes \Mod_A\simeq \LMod_{B\otimes A}(\Mod_{A})$ containing all objects $F$ such that $\Hom_{B\otimes A}(F,G)$ is perfect. To prove our statement, it is enough to show that $\Perf(B\otimes A)$ is contained in $\m{W}$. For this, it suffices to show that $B \otimes A$ belongs to $\m{W}$.
    
The functor $\Hom_{B \otimes A}(B \otimes A, -)$ is precisely the forgetful functor (i.e., the right adjoint to $B \otimes A \otimes_{A} -$). Since $B$ is smooth and proper, the forgetful functor preserves and reflects dualizable objects (see \cite[4.6.4.4 and 4.6.4.12]{HA}), which in this setting coincide with perfect complexes over $B\otimes A$.
\end{proof}

\begin{remark}
We can identify the \cat $\mathscr{M}_{B}(A)$ with the \cat of perfect $A$-modules equipped with an action of $B$, which we denote by $\LMod_{B}(\Perf(A))$; this is, in turn, equivalent to $\LMod_{B\otimes A}(\Perf(A))$. When $B$ is smooth and proper, $B\otimes A$ is smooth and proper over $A$. Therefore, in this case, by \Cref{lem: perfects on a smooth and proper algebra}, the \cat $\LMod_{B}(\Perf(A))$ coincides with $\Perf(B\otimes A)$.
\end{remark}

In this setting, \Cref{def: simple} translates to the following:

\begin{definition}\label{def:simple algebraic setting}
    We say that $F \in \LMod_{B\otimes A}$ is \textit{simple} if it is a perfect complex and the canonical map 
    \begin{equation*}
        \alpha \colon A \to \Hom_{B \otimes A}(F,F)
    \end{equation*}
    constructed as in \Cref{rmk:map simple}, is such that its cofiber $\cofib(\alpha)$ has Tor-amplitude strictly less than 0. 
\end{definition}

\begin{remark}\label{rmk: Hom and tensor product}
     Let $F$ and $G$ be perfect complexes in $\LMod_{B \otimes A}$, and let $A'$ be in $\Mod_{A}$. Then
     \begin{equation}\label{eq: Hom and tensor product}
    \Hom_{B \otimes A}(F,G) \otimes_{A} A' \simeq \Hom_{B \otimes A'}(F \otimes_{A} A', G \otimes_{A} A').
    \end{equation}
    By adjunction, providing the map in \eqref{eq: Hom and tensor product} is equivalent to providing a map:
    \begin{equation*}
        (\Hom_{B \otimes A}(F,G) \otimes_{A} A' \otimes_{A'} F) \otimes_{A} A' \to G \otimes_{A} A'
    \end{equation*}
    which we can take to be the counit. To prove that the map in \eqref{eq: Hom and tensor product} is an equivalence, we consider the smallest full stable subcategory of $\LMod_{B \otimes A}$ containing modules $F$ for which \eqref{eq: Hom and tensor product} is an equivalence. Since this category is stable under finite limits, finite colimits, and retracts, it must coincide with the category of perfect objects in $\LMod_{B \otimes A}$ (see \cite[Definition 7.2.4.1]{HA}). Thus, \eqref{eq: Hom and tensor product} is an equivalence for every perfect complex. 
\end{remark}

\begin{lemma}\label{lem: caracterization of simple algebraic setting}
     The following statements are equivalent: 
     \begin{enumerate}
         \item $F \in \Perf(B \otimes A)$ is simple.
         \item For every morphism $A \to A'$ of derived commutative algebras, $F\otimes_{A}A'$ is simple; i.e., the cofiber of the canonical map induced by taking the tensor product in \Cref{def:simple algebraic setting},
         \begin{equation*}
             \cofib(A' \to \Hom_{B \otimes A'}(F \otimes_{A} A', F \otimes_{A} A'))
         \end{equation*}
        has Tor-amplitude strictly less than 0.
        \item The same statement holds as above, but with $A'$ underived. 
     \end{enumerate}
\end{lemma}

\begin{proof}
  $(1) \Rightarrow (2)$: This follows from the fact that Tor-amplitude is stable under base change. Indeed, by \Cref{rmk: Hom and tensor product}, we have:
  \begin{align*}
      \cofib(A \to \Hom_{B \otimes A}(F,F)) \otimes_{A} A' & \simeq \cofib(A \otimes_{A} A' \to \Hom_{B \otimes A}(F,F) \otimes_{A} A' )\\
      & \simeq \cofib(A' \to \Hom_{B \otimes A'}(F \otimes_{A} A', F \otimes_{A} A')).
  \end{align*}
  The converse immediately follows by taking $A' = A$ and the identity map. \\
  $(3) \Rightarrow (2)$: We observe that $M \in \Mod_{A'}$ has Tor-amplitude strictly less than 0 if and only if this property holds for $M \otimes_{A'} \pi_{0}A'$. We apply this to $M = \cofib(A' \to \Hom_{B \otimes A'}(F \otimes_{A} A', F \otimes_{A} A'))$. The opposite direction is trivial.
\end{proof}

\begin{lemma}\label{lem: simple perfect complexes}
    Let $F \in \Perf(B \otimes A)$. Then $F$ is simple if and only if for every map $A \to A'$ of commutative algebras with $A'$ underived, we have:
    \begin{equation*}
        \pi_{i}\Hom_{B \otimes A'}(F \otimes_{A} A', F \otimes_{A} A') = 0 \quad \forall i \geq 1
    \end{equation*}
    and
    \begin{equation*}
        \pi_{0}\Hom_{B \otimes A'}(F \otimes_{A} A', F \otimes_{A} A') \simeq A'.
    \end{equation*}
\end{lemma}

\begin{proof}
    The forward implication ($\Rightarrow$) follows from statement (3) in \Cref{lem: caracterization of simple algebraic setting}, via the exact sequence induced by the fiber sequence: 
    \begin{equation*}
        A' \to \Hom_{B \otimes A'}(F \otimes_{A} A', F \otimes_{A} A') \to \cofib.
    \end{equation*}
    For the converse implication ($\Leftarrow$), by statement (3) in \Cref{lem: caracterization of simple algebraic setting}, we may assume that $A$ is underived. By hypothesis, the cofiber has homological amplitude concentrated in degrees strictly less than 0. Consequently, for every map $A \to k$, where $k$ is a field, the cofiber 
    \begin{equation*}
        \cofib(k \to \Hom_{B \otimes k}(F \otimes_{A} k, F \otimes_{A} k))
    \end{equation*}
    has homological amplitude in degrees strictly less than 0 (which, over a field, coincides with Tor-amplitude). We then conclude by applying the local criterion of flatness (see \cite[Corollary 6.1.4.7]{SAG}) to $M = \cofib(A \to \Hom_{B \otimes A}(F,F))$. 
\end{proof}

\begin{remark}\label{rmk: semplici sui punti}
    Note that the same proof works if we only consider maps from $A$ to the residue fields of $\pi_{0}A$.
\end{remark}

As an immediate corollary, we deduce the following:

\begin{corollary}\label{cor: simple objects algebraic underived}
    If $A$ is an underived commutative algebra, then $F \in \Perf(B\otimes A)$ is simple if and only if the canonical map 
    \[
    A \to \tau_{\geq 0}\Hom_{B \otimes A}(F,F)
    \]
    is an equivalence.
\end{corollary}

This corollary establishes the equivalence between \cite[Definition 4.1]{tova} and our \Cref{def:simple algebraic setting}. We also emphasize that when $A$ is a derived commutative algebra, the characterization presented in the preceding corollary is not appropriate. A naive definition of a simple perfect complex might stipulate that the canonical map 
\[
A \to \tau_{\geq 0}\Hom_{B \otimes A}(F,F)
\]
is an equivalence. However, when $A$ is derived, this condition is not stable under base change, and \Cref{lem: simple perfect complexes} fails, as demonstrated by the following example.

\begin{example}
    Let us consider $\Mod_{k}$ equipped with an action of $k[T]$ induced by the natural map $k[T] \to k$. We define $k[\epsilon]$ to be $k \oplus k[1]$, and we consider $\Mod_{k[\epsilon]}$ equipped with an action of $k[T, \epsilon] = k[T] \oplus k[T][1]$. The inclusion map induces a base-change functor $f^* \colon \Mod_{k} \to \Mod_{k[\epsilon]}$ given by $- \otimes_{k} k[\epsilon]$. Note that $k \in \Mod_{k}$ is simple with respect to the action of $k[T]$, since $\Hom_{k[T]}(k,k) \simeq k \oplus k[-1]$. However, $f^*(k) = k[\epsilon]$ fails to be simple under the naive definition because $\Hom_{k[T, \epsilon]}(f^*k, f^*k) \simeq \Hom_{k[T]}(k,k) \otimes_{k} k[\epsilon] \simeq k[1] \oplus k^2 \oplus k[-1]$. In particular, the degree $0$ part does not reduce merely to $k$.
\end{example}

\subsection{Simple perfect complexes over a rigid analytic variety}\label{sec2.3}\leavevmode \\

We define the sheaf of small \cates
\begin{equation}\label{eq: perf as sheaf of cat}
    \Perf(-) \colon \Dafd_{k} \to \mathrm{Cat}_{\infty}
\end{equation}
to be the functor that sends a derived affinoid $S$ to the \cat of perfect complexes over $S$, denoted by $\Perf(S)$ (see \cite[Definition 7.5]{derhom}). We denote by $\uline{\Perf}$ the associated derived analytic stack. 

We begin by recalling \cite[Theorem 3.4]{derhom}. We fix a derived affinoid space $S$ with $A$ its algebra of global sections and we consider the following diagram.
\begin{equation*}
\begin{tikzcd}
    \Mod_{A} \arrow[rr, "\Psi"]                   &  & \Mod_{\mathscr{O}_{S}}                   \\
    \Perf(A) \arrow[rr] \arrow[u, hook] &  & \Perf(S). \arrow[u, hook]
\end{tikzcd}
\end{equation*}
Here the functor $\Psi$ sends an $A$-module $M$ to the sheaf $\Psi(M)$, which, when evaluated on an \'etale affinoid $\hat{S} \to S$ with $B\coloneqq \Gamma(\hat{S}, \m{O}_{\hat{S}})$, yields the module $M \widehat{\otimes}_{A} B$, where $\widehat{\otimes}$ denotes the completed tensor product. The functor $\Psi$ is symmetric monoidal and its inverse is given by the global sections functor.
Notably, it follows from \cite[Theorem 3.4 and Lemma 7.6]{derhom} that the functor $\Psi$ restricts to perfect complexes; in this case, the completion is unnecessary as these modules are finitely generated. Using the results of \cite{derhom}, we moreover deduce the following:

\begin{theorem}[{\cite[Theorem 3.4]{derhom}, Tate acyclicity}]\label{thm: Tate acyclicity}
    The functor $\Psi$ is t-exact and induces a symmetric monoidal equivalence between the \cates of perfect complexes, with its inverse given by the global sections functor $\Gamma(-,S)$.
\end{theorem}
\begin{proof}
    We can deduce that $\Psi$ is an equivalence using \cite[Theorem 3.4 and Lemma 7.6]{derhom}. To prove that $\Psi$ is symmetric monoidal, we observe that in the proof of \cite[Theorem 3.4]{derhom}, for $M \in \Perf(A)$, $\Psi(M)$ is identified with the presheaf 
    \[
    \mathscr{O}_{S}(-)\otimes M \colon \mathrm{dAfd}_{/S} \to \Mod_{\mathscr{O}_{S}},
    \]
    which is, moreover, a hypersheaf. This implies that $\Psi(-)$ is symmetric monoidal. 
\end{proof}

\begin{proposition}[{\cite{mathew}}]\label{prop: descent for perf}
    The functor in \eqref{eq: perf as sheaf of cat} satisfies descent for the flat topology.
\end{proposition}
\begin{proof}
    This follows from \cite[Theorem 7.8]{mathew} and \Cref{thm: Tate acyclicity}.  
\end{proof}

\begin{remark}
    Note that the functor $\Perf(-)$ is also a sheaf of stable, idempotent-complete \cates. This follows from \cite[Proposition 5.1.4.9]{HTT} and \cite[Theorem 1.1.4.4]{HA}.
\end{remark}

\begin{definition}
    We define the analytic stack of perfect complexes over a derived analytic stack $X$, denoted by $\uline{\Perf}_{X}$, as:
    \[
    \uline{\Perf}_{X} \coloneqq \uline{\mathrm{Map}}(X, \uline{\Perf}),
    \]
    where $\uline{\mathrm{Map}}$ denotes the analytic mapping stack as defined in \cite{derhom}.
\end{definition}

\begin{remark}
    We note that the functor of points of $\uline{\Perf}_{X}$ can be described explicitly by sending an affinoid space $S$ to the \cat $\Perf(X \times S)^{\simeq}$.
\end{remark}

We will now formulate the notion of a simple object in $\uline{\Perf}_{X}(S)$, where $X$ is a smooth and proper rigid analytic variety and $S$ is a derived affinoid with $A=\Gamma(S,\m{O}_{S})$.\\ 
To do this, we apply \Cref{rmk:map simple} with $\m{V}=\Mod_{\m{O}_{S}}$ and $\m{C}=\Mod_{\m{O}_{X\times S}}$, where the action is induced by sending the pair $(\m{F}\in \Mod_{\m{O}_{S}}, \m{G}\in \Mod_{\m{O}_{X \times S}})$ to $p^{*}\m{F}\otimes_{X\times S} \m{G}$. Here, $p \colon X \times S \to S$ denotes the projection. 
In this case, the adjunction can be written more explicitly. Indeed, by composing the tensor-Hom adjunction for $\m{O}_{X \times S}$-modules with the adjunction between $p^{*}$ and $p_{*}$, we obtain the following adjunction:

\begin{equation*}
    \begin{tikzcd}
	{\Mod_{\m{O}_{X\times S}}} &&& {\Mod_{\m{O}_{S}}}
	\arrow[""{name=0, anchor=center, inner sep=0}, "{p^{*}(-)\otimes_{X\times S}\m{G}}", shift left=2, from=1-4, to=1-1]
	\arrow[""{name=1, anchor=center, inner sep=0}, "{p_{*}\m{Hom}_{X\times S}(\m{G},-)}", shift left=2, from=1-1, to=1-4]
	\arrow["\dashv"{anchor=center, rotate=+90}, draw=none, from=0, to=1]
\end{tikzcd}
\end{equation*}
where we denote by $\mathscr{Hom}(-,-)$ the internal Hom sheaf in $\Mod_{\mathscr{O}_{S}}$.

\begin{lemma}\label{lemma: hom perfect}
    Let $Y$ be a rigid analytic variety and let $\m{F}, \m{G}$ be perfect complexes over $Y$. Then $\mathscr{Hom}_{\m{O}_{Y}}(\m{F},\m{G})$ is also a perfect complex.
\end{lemma}
\begin{proof}
    The statement is local on $Y$. Locally, we may assume that $\m{F}$ is a retract of a finite colimit of copies of $\m{O}_{Y}$. This implies that, locally, $\m{Hom}_{\m{O}_{Y}}(\m{F},\m{G})$ is a retract of a finite colimit of copies of $\m{Hom}_{\m{O}_{Y}}(\m{O}_{Y},\m{G})\simeq \m{G}$, which is a perfect complex. Therefore, $\m{Hom}_{X\times S}(\m{F},\m{G})$ is also a perfect complex.
\end{proof}

\begin{corollary}\label{cor: right adjoint for sheaves}
    Let $\m{F} \in \Perf(X \times S)$ and let $p_{S} \colon X \times S \to S$ be the projection map. Then there is an adjunction: 
    \begin{equation*}
       \begin{tikzcd}
	{\Perf(X\times S)} &&& {\Perf(S).}
	\arrow[""{name=0, anchor=center, inner sep=0}, "{\m{F}\otimes_{X\times S}p^{\ast}_{S}(-)}", shift left=2, from=1-4, to=1-1]
	\arrow[""{name=1, anchor=center, inner sep=0}, "{{p_{S\ast}\m{Hom}_{X\times S}(\m{F},-)}}", shift left=2, from=1-1, to=1-4]
	\arrow["\dashv"{anchor=center, rotate=+90}, draw=none, from=0, to=1]
    \end{tikzcd}
    \end{equation*}
\end{corollary}
\begin{proof}
    Starting from the adjunction between $\Mod_{\m{O}_{X\times S}}$ and $\Mod_{\m{O}_{S}}$, we observe that since $X$ is smooth and proper, the functors $p^{\ast}_{S}$ and $p_{S\ast}$ preserve perfect objects \cite[Proposition 7.8]{derhom}. Furthermore, the tensor product of two perfect complexes is again a perfect complex. Thus, we only need to prove that $\m{H}om_{X\times S}(\m{F},\m{G})$ is a perfect complex whenever $\m{G}$ is perfect. This follows immediately from \Cref{lemma: hom perfect}.
\end{proof}

We denote the right adjoint in \Cref{cor: right adjoint for sheaves} by $\uline{\m{Hom}}_{X\times S}(\m{F},-)$ and use it to describe simple objects as in \Cref{def: simple}. In particular, we have the following definition:

\begin{definition}\label{def: simple analytic}
    An object $\m{F}\in \Perf(X\times S)$ is \emph{simple} if the cofiber of the canonical map (as constructed in \Cref{rmk:map simple}),
    \begin{equation*}
        \m{O}_{S} \to \uline{\m{H}om}_{X\times S}(\m{F},\m{F}),
    \end{equation*}
   has Tor-amplitude in strictly negative degrees. 
\end{definition}

\begin{remark}\label{rmk: diagramma pullback sul punto}
    Let $S$ be a derived affinoid and $s$ a closed point in $S$ (i.e., a closed point in the associated underived affinoid $\pi_{0}S$) classified by a map $s \colon \Sp(K) \to S$. Consider the following Cartesian diagram: 
    \begin{equation*}
        \begin{tikzcd}
X\times \Sp(K) \arrow[rr, "i_{s} \coloneqq id \times s"] \arrow[dd, "q"'] &  & X \times S \arrow[dd, "p"] \\
                                               &  &                            \\
{\Sp(K)} \arrow[rr, "s"']                      &  & S                         
\end{tikzcd}
    \end{equation*}
    where $p$ and $q$ are the projections. We have that 
    \begin{equation*}
        s^{*}\uline{\m{H}om}_{X \times S}(\m{F},\m{F}) \coloneqq s^{*}p_{*}\m{H}om_{X\times S}(\m{F},\m{F})\simeq q_{*}i_{s}^{*}\m{H}om_{X\times S}(\m{F},\m{F}).
    \end{equation*}
    Now, since a perfect complex is dualizable and the functor $i_{s}^{*}$ is symmetric monoidal, we obtain:
    \begin{equation*}
        s^{*}\uline{\m{Hom}}_{X \times S}(\m{F},\m{F})\simeq q_{\ast}\m{Hom}_{X\times \Sp(K)}(i_{s}^{*}\m{F},i_{s}^{*}\m{F}) \coloneqq \uline{\m{Hom}}_{X\times \Sp(K)}(i_{s}^{*}\m{F},i_{s}^{*}\m{F}).
    \end{equation*}
    In particular, taking global sections of the map in \Cref{def: simple analytic}, we obtain a map:
    \begin{equation}\label{eq: canonical map dal punto}
        K \to \Hom_{X\times \Sp(K)}(i_{s}^{\ast}\m{F},i_{s}^{\ast}\m{F}),
    \end{equation}
    where $K$ is the finite field extension of $k$ corresponding to the point $s$.
\end{remark}

\begin{lemma}\label{lem: caratterization simple analytic}
    Let $\m{F} \in \Perf(X \times S)$. Using the notation from \Cref{rmk: diagramma pullback sul punto}, $\m{F}$ is simple if and only if for every closed point $s\colon \Sp(K) \to S$, the sheaf $i_{s}^{\ast}\m{F}$ is simple (i.e., the canonical map
    \begin{equation}\label{eq: map simple relative to the point}
        \mathscr{O}_{\Sp(K)} \to q_{\ast}\m{Hom}_{X\times \Sp(K)}(i_{s}^{*}\m{F},i_{s}^{*}\m{F}) \coloneqq \uline{\m{Hom}}_{X\times \Sp(K)}(i_{s}^{*}\m{F},i_{s}^{*}\m{F})
    \end{equation}
    has its cofiber in strictly negative homological degrees).
    This is equivalent to requiring that the map constructed in \eqref{eq: canonical map dal punto},  
    \begin{equation*}
        K \to \Hom_{X\times \Sp(K)}(i_{s}^{\ast}\m{F},i_{s}^{\ast}\m{F}),
    \end{equation*}
    has its cofiber in strictly negative homological degrees.
\end{lemma}
\begin{proof}
    We first observe that applying \Cref{thm: Tate acyclicity} to the cofiber of the map in \Cref{def: simple analytic} yields an object in $\Perf(A)$ given by
    \[
    \mathrm{cofib}(A \to \mathrm{Hom}_{X\times S}(\m{F},\m{F})).
    \]
    By \Cref{thm: Tate acyclicity}, we deduce that this object has Tor-amplitude in strictly negative degrees if and only if $\m{F}$ is simple. Consequently, for every closed point 
    \[
    s\colon \Sp(K) \to S,
    \]
    the cofiber of the map in \eqref{eq: map simple relative to the point} has Tor-amplitude in strictly negative degrees if and only if the same holds for the map in \eqref{eq: canonical map dal punto}. Therefore, showing that $\m{F}$ is simple if and only if $i^{\ast}_{s}\m{F}$ is simple reduces to \Cref{lem: caracterization of simple algebraic setting}. Indeed, upon taking global sections (which preserves simple objects), the pullback $s^{\ast}$ can be identified with the functor $-\otimes_{A}K$. This yields exactly the map in \eqref{eq: canonical map dal punto}. 
\end{proof}

In this way, we can see that our \Cref{def: simple analytic} is exactly the non-Archimedean analogue of \cite[Definition 5.2]{tova}.

\begin{corollary}\label{cor: simple analytic underived}
    If $S$ is underived, then $\m{F}$ is simple if and only if the map from \Cref{def: simple analytic} induces an equivalence: 
    \begin{equation*}
        \m{O}_{S} \simeq \tau_{\geq 0}\uline{\m{H}om}_{X\times S}(\m{F},\m{F}).
    \end{equation*}
\end{corollary}
\begin{proof}
    If $S$ is underived, then $\m{O}_{S}$ is concentrated in degree $0$. The result now follows from \Cref{lem: caratterization simple analytic} and \Cref{lem: simple perfect complexes}.
\end{proof}

\subsection{The rigid analytic stack of objects in an \texorpdfstring{$\infty$-category}{\cat}}\label{sec2.4}
\leavevmode \\

We begin this section by defining the (rigid) derived analytic stack of objects in an $\infty$-category.

\begin{definition}
    The functor of points $\m{AnM}_{\m{C}} \colon \Dafd_{k}^{\mathrm{op}} \to \m{S}$ is defined by sending a derived affinoid $S$ to the space of functors $\Fun^{\mathrm{st}}_{k}(\m{C}^{\omega}, \Perf(S))^{\simeq}$.
\end{definition}

As a corollary of \Cref{prop: descent for perf}, we obtain the following result:
\begin{corollary}
     The functor $\m{AnM}_{\m{C}}$ is a derived analytic stack.
\end{corollary}

We will now define the substack of simple objects.\\
Consider a derived affinoid space $S$ and let $A \coloneqq \Gamma(S, \m{O}_{S})$. To define simple objects in $\mathscr{AnM}_{\Mod_{B}}(S)$, we will utilize the definition of simple objects established above for the algebraic case, combined with Tate acyclicity, as follows.

\begin{remark}\label{rmk: AnM(B) and LModB}
When $\m{C}$ is equivalent to $\LMod_{B}$ for a smooth and proper $\mathbb{E}_{1}$-algebra $B$, we have
\[
\m{AnM}_{\m{C}}(S) \coloneqq \Fun^{\mathrm{st}}_{k}(\Perf(B), \Perf(S)),
\]
which classifies perfect $\m{O}_{S}$-modules equipped with an action of $B$. For this reason, we can identify this \cat with $\LMod_{B}(\Perf(S))$.
\end{remark}

\begin{corollary}\label{cor: tate acyclicity for relative to an algebra}
    The functor $\Psi(-)$ from \Cref{thm: Tate acyclicity} induces an equivalence between $\LMod_{B}(\Perf(A))$ and $\LMod_{B}(\Perf(S))$.
\end{corollary}
\begin{proof}
    This is a consequence of \Cref{thm: Tate acyclicity} and \Cref{lem:symmetric_monoidal_functors_and_modules}, taking into account \Cref{rmk: AnM(B) and LModB} and \Cref{rmk: MB vs LModB}.
\end{proof}

We observe that since $\Psi(-)$ is symmetric monoidal, by using \Cref{lem:symmetric_monoidal_functors_and_modules} we also deduce that the \cat $\m{AnM}_{\m{C}}(S)$ is equivalent to $\LMod_{\Psi(B\otimes A)}(\Perf(S))$. 

\begin{proposition}
    For every object $\mathscr{F}$ in $\LMod_{B}(\Perf(S))$, the functor 
    \[
    \mathscr{F} \otimes_{\mathscr{O}_{S}} - \colon \Perf(S) \to \LMod_{B}(\Perf(S))
    \]
    admits a right adjoint.
\end{proposition}
\begin{proof}
    Using \Cref{thm: Tate acyclicity}, \Cref{cor: tate acyclicity for relative to an algebra}, and \Cref{prop:afjunction for simple objects algebra}, we obtain the following commutative diagram:
    \begin{equation}\label{eq: adjunction affinoid setting}
    \begin{tikzcd}
	{\LMod_{B}(\Perf(S))} && {\LMod_{B}(\Perf(A))} \\
	\\
	{\Perf(S)} && {\Perf(A).}
	\arrow["{\simeq }"', from=1-3, to=1-1]
	\arrow[""{name=0, anchor=center, inner sep=0}, "{\Hom_{B\otimes A}(\Gamma(\mathscr{F},S),-)}", shift left=3, from=1-3, to=3-3]
	\arrow["{\mathscr{F}\otimes-}", from=3-1, to=1-1]
	\arrow[""{name=1, anchor=center, inner sep=0}, "{\Gamma(\mathscr{F},S) \otimes -}", shift left=3, from=3-3, to=1-3]
	\arrow["{\simeq }", from=3-3, to=3-1]
	\arrow["\dashv"{anchor=center}, draw=none, from=1, to=0]
\end{tikzcd}  
\end{equation}
 From \eqref{eq: adjunction affinoid setting}, we deduce the existence of the required right adjoint.
\end{proof}

\begin{definition}\label{def: non commutative hom affine}
    For $\m{F} \in \Mod_{B}(\Perf(S))$, we denote by $\uline{\Hom}_{B}(\m{F},-)$ the right adjoint to $\m{F} \otimes_{\m{O}_{S}} -$.
\end{definition}

\begin{definition}\label{def: simple affine setting}
    We say that $\m{F} \in \LMod_{B}(\Perf(S))$ is \textit{simple} if the cofiber of the canonical map
    \begin{equation*}
       \m{O}_{S} \to \uline{\Hom}_{B}(\m{F},\m{F})
    \end{equation*}
    has Tor-amplitude concentrated in homological degrees strictly less than 0. 
\end{definition}

\begin{proposition}\label{rmk: explicit definition simple relative to an algebra}
    Let $\m{F}$ be an object in $\LMod_{B}(\Perf(S))$, and let
    \[
    F \coloneqq \Gamma(\m{F},S)
    \]
    be the corresponding $A$-module of global sections. 
    Then $F$ is simple in the sense of \Cref{def:simple algebraic setting} if and only if $\m{F}$ is simple in the sense of \Cref{def: simple affine setting}.
\end{proposition}
\begin{proof}
    By \Cref{thm: Tate acyclicity}, we can assume without loss of generality that $\m{F} \simeq \Psi(F)$.
    Consider the canonical map from \Cref{def:simple algebraic setting}:
        \begin{equation*}
            \begin{tikzcd}
A \arrow[r, "\alpha "] & {\Hom_{B\otimes A}(F,F)}
\end{tikzcd}
        \end{equation*}
        This induces the following diagram in $\Mod_{\mathscr{O}_{S}}$, where $\beta$ is the map introduced in \Cref{def: simple affine setting}:
        \begin{equation*}
        \begin{tikzcd}
\Psi(A) \arrow[d, "\simeq "'] \arrow[rr, "\Psi(\alpha)"] &  & {\Psi(\Hom_{B\otimes A}(F,F))} \arrow[d, "\simeq"]             \\
\mathscr{O}_{S} \arrow[rr, "\beta"']                   &  & {\uline{\Hom}_{B}(\m{F},\m{F})}
\end{tikzcd}
        \end{equation*}
The vertical map on the left is an equivalence (see \Cref{def: simple affine setting} and \eqref{eq: adjunction affinoid setting}).
From this, we deduce that $\Psi(\alpha)$ is an equivalence if and only if $\beta$ is an equivalence.
Using \Cref{thm: Tate acyclicity}, we conclude that $\alpha$ is an equivalence if and only if $\beta$ is an equivalence. 
\end{proof}

\begin{lemma}\label{lem: simple analitic relative to B stable under BC}
    Let $\m{F}$ be an object in $\LMod_{B}(\Perf(S))$. Then, for every map of derived affinoids 
    \[
    i\colon S' \to S,
    \]
    $\m{F}$ is simple if and only if $i^{\ast}\m{F}$ is simple.
\end{lemma}
\begin{proof}
    In light of \Cref{rmk: explicit definition simple relative to an algebra}, this lemma can be deduced directly from \Cref{lem: caracterization of simple algebraic setting}. 
\end{proof}

\begin{lemma}\label{lem: caracterization simple affine}
    Let $\m{F} \in \mathscr{AnM}_{B}(S)$, where $B$ is a smooth and proper algebra. Then $\m{F}$ is simple if and only if for every closed point $i \colon s \to S$ associated to a field $K$, the cofiber of the canonical map 
    \begin{equation*}
        K \to \Hom_{B \otimes K}(i^{*}\m{F}, i^{*}\m{F})
    \end{equation*}
    (as in \Cref{rmk:map simple}) induces the following equivalences: 
    \begin{equation*}
        \pi_{0}\Hom_{B \otimes K}(i^{*}\m{F}, i^{*}\m{F}) \simeq K,
    \end{equation*}
    \begin{equation*}
        \pi_{j}\Hom_{B \otimes K}(i^{*}\m{F}, i^{*}\m{F}) \simeq 0 \quad \text{for every } j > 0.
    \end{equation*}
\end{lemma}

\begin{proof}
    By passing to global sections and using \Cref{rmk: explicit definition simple relative to an algebra}, we can reduce to the algebraic setting. Our statement then corresponds to \Cref{lem: simple perfect complexes} applied to every map $A \to K$, where $K$ is a finite field extension (see also \Cref{rmk: semplici sui punti}).
\end{proof}

This proposition demonstrates that \Cref{def: simple affine setting} is the non-Archimedean derived analogue of \cite[Definition 4.9]{tova}.

\begin{corollary}\label{cor: caracterization simple affinoid underived}
    If $S$ is an underived affinoid, then $\m{F} \in \Mod_{B}(\Perf(S))$ is simple if and only if the canonical map from \Cref{def: simple affine setting} induces an equivalence: 
    \begin{equation*}
        \m{O}_{S} \simeq \tau_{\geq 0} \uline{\Hom}_{B}(\m{F}, \m{F}).
    \end{equation*}
\end{corollary}

\begin{proof}
   If $S$ is underived, then $\m{O}_{S}$ is concentrated in degree 0. Our statement now follows from \Cref{lem: caracterization simple affine}.
\end{proof}

\subsection{Comparison between the stacks \texorpdfstring{$\m{M}_{B}$}{MB} and \texorpdfstring{$\m{AnM}_{B}$}{AnMB}}\label{secComparison1}\leavevmode \\

The aim of this subsection is to prove that if $\mathscr{C}$ is a smooth and proper \cat, then there is an equivalence between $\m{AnM}_{\mathscr{C}}$ and $(\m{M}_{\mathscr{C}})^{\mathrm{an}}$, and that this equivalence moreover preserves simple objects (see \Cref{prop: MB sim vs AnMB sim}). The proof is based on comparing the functors of points of the two stacks. We begin by introducing the main tool we need (i.e., \cite[Theorem 3.10]{anmap}).

Let $Y$ be a derived algebraic stack over a non-Archimedean field $k$, and let $Y^{\mathrm{an}}$ be its analytification as in \cite{anmap}. We consider the following functor: 

\begin{definition}
We denote by $G(Y)$ the functor:
\begin{equation*}
\begin{tikzcd}[row sep=0.01in]
G_\ast(Y) \colon \Dafd \arrow[rr] &  & \mathscr{S} \\
S \arrow[rr]            &  & Y(\Spec(\Gamma(S,\m{O}_{S})))   
\end{tikzcd}
\end{equation*}
In the notation of \cite{anmap}, this presheaf is denoted by $G_{\ast}^{p}(Y)$, and its sheafification is denoted by $G_{\ast}^{s}(Y)$, where $\ast$ is the final object in $\mathrm{dAnSt}_{k}$.
\end{definition}

\begin{theorem}[{\cite[Theorem 3.10]{anmap}}]\label{thm: teorema mauro} 
    Let $Y$ be a locally geometric derived stack locally almost of finite presentation. Then there is an equivalence of presheaves:
    \begin{equation}
        G_{\ast}(Y) \simeq Y^{\mathrm{an}}.
    \end{equation}
\end{theorem}

Using \Cref{thm: teorema mauro}, it is straightforward to see that the functor of points of $(\mathscr{M}_{\mathscr{C}})^{\mathrm{an}}$ coincides with that of $\mathscr{AnM}_{\mathscr{C}}$. We will first construct a map between these two derived stacks, and subsequently prove that it is an equivalence.

\begin{construction} 
Let $\m{C}$ be a smooth and proper \cat in $\mathrm{Pr}^{\mathrm{L},\omega}$. We consider the composition: 
\begin{equation}\label{eq: analtification moduli}
    \begin{tikzcd}
	{\m{AnM}_{\mathscr{C}}} && {G(\m{M}_{\mathscr{C}})} && { (\m{M}_{\mathscr{C}})^{\mathrm{an}}}
	\arrow[from=1-1, to=1-3]
	\arrow["\alpha"', from=1-3, to=1-5]
\end{tikzcd}
\end{equation}
where $\alpha$ is the map considered in \cite[Theorem 3.10]{anmap}, and the first map is induced by the global sections functor, as in \Cref{thm: Tate acyclicity}.
\end{construction}

\begin{proposition}\label{cor: analytification and analytic stack}
    If $\m{C}$ is a smooth and proper \cat in $\mathrm{Pr}^{\mathrm{L},\omega}_{\mathrm{st}}$, then the map in \eqref{eq: analtification moduli} is an equivalence.
\end{proposition}
\begin{proof}
We will prove that both maps 
\[
(\m{M}_{\mathscr{C}})^{\mathrm{an}} \to G(\m{M}_{\mathscr{C}})
\]
and 
\[
 \m{AnM}_{\mathscr{C}} \to G(\m{M}_{\mathscr{C}})
\]
are equivalences.
To prove that the first map is an equivalence, we observe that \cite[Corollary 3.15]{modg} (since $\mathscr{C}$ is smooth and proper) states that the functor $\mathscr{M}_{\mathscr{C}}$ is a derived geometric stack locally of finite presentation; therefore, we can apply \Cref{thm: teorema mauro}. Now we move to the second equivalence, where can then easily compare the functors of points of $G(\mathscr{M}_{\mathscr{C}})$ and $\mathscr{AnM}_{\mathscr{C}}$. Indeed, they are equivalent by Tate acyclicity (\Cref{thm: Tate acyclicity}).
\end{proof}

Note that this result allows us to describe the functor of points of the analytification of $\mathscr{M}_{\mathscr{C}}$, which is otherwise only defined as a left Kan extension. Indeed, we now have an explicit description for the functor of points via the stack $\mathscr{AnM}_{\mathscr{C}}$.
Since we always assume $\mathscr{C}$ to be a smooth and proper \cat in $\mathrm{Pr}^{\mathrm{L},\omega}_{\mathrm{st}}$, by \Cref{cor: categories vs algebras} we may assume without loss of generality that $\mathscr{C}$ is of the form $\LMod_{B}$, where $B$ is a smooth and proper $\mathbb{E}_{1}$-algebra. We now aim to prove that the equivalence in \Cref{cor: analytification and analytic stack} restricts to simple objects.

\begin{definition}\label{def: stack of simple objects algebraic setting}
    We define $\mathscr{M}_{B}^{\mathrm{sim}}$ as the substack of $\mathscr{M}_{B}$ that sends a derived commutative algebra $A$ to the groupoid of simple objects $\m{M}_{B}^{\mathrm{sim}}(A)^{\simeq}$, as defined in \Cref{def:simple algebraic setting}.  
\end{definition}

Using \Cref{lem: caracterization of simple algebraic setting}, we can readily see that $\mathscr{M}_{B}^{\mathrm{sim}}$ is a substack of $\mathscr{M}_{B}$.

\begin{remark}\label{rmk: artin 1-stack of simple objects}
Using \Cref{cor: simple objects algebraic underived}, we deduce that the functor $\mathscr{M}^{\mathrm{sim}}_{B}$, when restricted to underived commutative algebras (denoted by $t_{0}\mathscr{M}^{\mathrm{sim}}_{B}$), is naturally a $1$-stack. Moreover, by \cite[Corollary 3.22]{modg}, it is an Artin $1$-stack. We also note that $t_0\mathscr{M}_{B}^{\mathrm{sim}}$ corresponds to the Artin $1$-stack of simple objects defined in \cite{modg}.
\end{remark}

\begin{definition}
    We define $\mathscr{AnM}_{B}^{\mathrm{sim}}$ as the functor that sends a derived affinoid space $S$ to the groupoid of simple objects $\m{AnM}_{B}^{\mathrm{sim}}(S)^{\simeq}$, as defined in \Cref{def: simple affine setting}.
\end{definition}

The functoriality of the above definition follows from \Cref{lem: simple analitic relative to B stable under BC}. Using \Cref{cor: caracterization simple affinoid underived}, we can see that the restriction of $\m{AnM}^{\mathrm{sim}}_{B}$ to underived affinoids, denoted by $t_0\m{AnM}_{B}^{\mathrm{sim}}$, is naturally a $1$-stack. 

\begin{proposition}\label{prop: MB sim vs AnMB sim}
    If $B$ is a smooth and proper algebra, then the equivalence of \Cref{cor: analytification and analytic stack} restricts to an equivalence of functors $(\mathscr{M}_{B}^{\mathrm{sim}})^{\mathrm{an}} \simeq \mathscr{AnM}_{B}^{\mathrm{sim}}$.
\end{proposition}
\begin{proof}
    We begin by observing that \Cref{cor: analytification and analytic stack}, applied with $\m{C}=\Mod_{B}$, yields the equivalence between $(\m{M}_{B})^{\mathrm{an}}$ and $\m{AnM}_{B}$. To establish the equivalence between the subfunctors of simple objects, we apply the same argument utilizing the functor $G$; we merely need to check that the functor $\Psi(-)$ preserves and reflects simple objects. This has been proven in \Cref{rmk: explicit definition simple relative to an algebra}. 
\end{proof}

In particular, this proposition implies that $\m{AnM}_{B}$ is also a locally geometric derived analytic stack.

\subsection{Comparison between \texorpdfstring{$\uline{\Perf}_{X}$}{PerfX} and \texorpdfstring{$\m{AnM}_{\Perf(X)}$}{AnMPerfX}}\label{secComparison2}

\begin{definition}\label{def: substack of perfect simple objects}
    We define the functor $\uline{\Perf}^{\mathrm{sim}}_{X}$ as the subfunctor of $\uline{\Perf}_{X}$ defined by sending a derived affinoid space $S$ to the groupoid of simple objects $\Perf(X\times S)^{\simeq}$.
\end{definition}

We observe that when $X$ is a proper (underived) rigid analytic variety, \Cref{def: substack of perfect simple objects} is well-posed since the pullback of a simple object remains simple. This follows from the base change property explained in \cite[Theorem 6.8]{derhom}.

The aim of this section is to prove the following proposition. 

\begin{proposition}\label{prop: Perf vs AnM}
    Let $X \to \Sp(k)$ be a proper, quasi-compact, and separated rigid analytic variety such that $\Perf(X)$ is a smooth and proper \cat over $\Perf(k)$. Then there is an equivalence 
    \[
    \uline{\Perf}_{X} \simeq \m{AnM}_{\Perf(X)}.
    \] 
    Furthermore, this equivalence restricts to an equivalence between simple objects (i.e., $\uline{\Perf}^{\, \mathrm{sim}}_{X} \simeq \m{AnM}_{\Perf(X)}^{\mathrm{sim}}$). 
\end{proposition}

To prove this proposition, we require several lemmas.

\begin{lemma}\label{lem: non commutative properness}
    Let $p\colon X \to \Sp(k)$ be a proper rigid analytic variety. If $\m{F}$ is a perfect complex, then $p_{*}(\m{F})$ is also a perfect complex in $\Mod_{k}$.
\end{lemma}
\begin{proof}
This lemma follows from \cite[Lemma 7.6]{derhom}. We merely need to show that $p$ has finite coherent cohomological dimension and finite Tor-amplitude. The morphism $p$ has finite coherent cohomological dimension because $X$ is of finite type. It also has finite Tor-amplitude because the target is a field and $X$ is underived. 
\end{proof}

\begin{lemma}\label{lem: compact generator for perf and adjuntion}
Let $p \colon X \to \Sp(K)$ be a proper rigid analytic variety over a non-Archimedean field $K$, such that $\Perf(X)$ is a smooth and proper, $K$-linear, idempotent-complete \cat. Then $\Perf(X)$ admits a strong generator $\m{E}$. Moreover, the pair of adjoint functors 
\begin{equation*}
    \begin{tikzcd}
	{\Perf(X)} &&& {\Perf(K)}
	\arrow[""{name=0, anchor=center, inner sep=0}, "{\Hom_{\m{O}_{X}}(\m{E},-)}"', shift right=2, from=1-1, to=1-4]
	\arrow[""{name=1, anchor=center, inner sep=0}, "{p^{\ast}(-)\otimes_{\m{O}_{X}}\m{E}}"', shift right=2, from=1-4, to=1-1]
	\arrow["\dashv"{anchor=center, rotate=-90}, draw=none, from=1, to=0]
\end{tikzcd}
\end{equation*}
induces an equivalence between $\Perf(X)$ and $\LMod_{B}(\Perf(K))$, where $B$ is the $\mathbb{E}_{1}$-algebra of endomorphisms of the generator $\m{E}$.
\end{lemma}
\begin{proof}
    We begin by observing that the adjunction is well-defined since $X$ is proper over $\Sp(K)$.
    By definition, $\mathrm{Ind}(\Perf(X))$ is a smooth and proper, presentable, and compactly generated \cat. The existence of a compact generator for $\mathrm{Ind}(\Perf(X))$ follows from \Cref{lem: compact generator}. The remainder of the proof relies on the Barr--Beck--Lurie Theorem \cite[Theorem 4.7.3.5]{HA}. Let $\m{E}$ be the compact generator of $\mathrm{Ind}(\Perf(X))$. The functor $\Hom_{\mathrm{Ind}(\Perf(X))}(\m{E},-)$ is conservative (since $\m{E}$ is a generator) and commutes with colimits (since $\m{E}$ is compact and the categories are stable). By the aforementioned theorem, the functor $\mathrm{Hom}_{\m{O}_{X}}(\m{E},-)$ induces an equivalence 
\[
\mathrm{Ind}(\Perf(X)) \simeq \LMod_{T}(\Mod_K),
\]  
where $T$ is the monad $\Hom_{\m{O}_{X}}(\m{E},\m{E}\otimes_{\m{O}_{X}}p^{\ast}(-))$. 
We will now identify the monad $T$ with $B\otimes_{K}-$.
    For every $K$-module $M$, there is a canonical map:
    \begin{equation}\label{eq: boh}
        \Hom_{X}(\m{E},\m{E})\otimes_{K}M \to \Hom_{X}(\m{E},\m{E}\otimes_{\m{O}_{X}} p^{\ast}M).
    \end{equation}
    Let $\m{W}$ be the full subcategory of $\Mod_K$ consisting of all objects for which \eqref{eq: boh} is an equivalence. We observe that $\m{W}$ contains $K$ and is closed under colimits, implying that it is the entire \cat $\Mod_{K}$. In this way, we obtain an equivalence: 
    \[
    \LMod_{T}(\Mod_K) \simeq \LMod_{B}.
    \]
    Since $B$ is a smooth and proper algebra in $\Mod_K$, we can identify the compact objects in $\LMod_{B}$ with $\LMod_{B}(\Perf(K))$. Because $X$ is proper, the functor $\Hom_{X}(\m{E},-)$ preserves perfect complexes, and thus it restricts to an equivalence between $\Perf(X)$ and $\LMod_{B}(\Perf(K))$. 
\end{proof}

\begin{lemma}\label{lem: compact generation after bc}
    Let $p\colon X \to \Sp(k)$ be a separated rigid analytic variety such that $\Perf(X)$ is smooth and proper over $\Perf(k)$. Let $K$ be a finite field extension of $k$. We denote by $X_K$ the fiber product rigid analytic variety in the following pullback diagram: 
    \[
    \begin{tikzcd}
	{X_K} && X \\
	\\
	{\Sp(K)} && {\Sp(k).}
	\arrow[from=1-1, to=1-3]
	\arrow[from=1-1, to=3-1]
	\arrow[from=1-3, to=3-3]
	\arrow[from=3-1, to=3-3]
\end{tikzcd}
    \]
    Then we have an equivalence:
    \[
    \Perf(X_K) \simeq \Perf(X)\otimes_{\Perf(k)}\Perf(K).
    \]
    Consequently, if $\Perf(X)$ is smooth and proper over $\Perf(k)$, then $\Perf(X_K)$ is smooth and proper over $\Perf(K)$.
\end{lemma}
\begin{proof}
    We begin by noting that $\Perf(K)$ is a dualizable \cat over $\Perf(k)$; by \cite[Theorem 3.4]{blumberg2013universal}, this is equivalent to saying that $\Perf(K)$ is smooth and proper over $\Perf(k)$. Indeed, $K$ is a finite separable field extension of $k$ (which is a field of characteristic 0), implying that the map 
    \[
    \Spec(K) \to \Spec(k)
    \]
    is \'etale. This ensures that $\Perf(K)$ is smooth over $\Perf(k)$. It is also proper because $K$ is a finitely generated projective $k$-module.
    When $X$ is isomorphic to an affinoid space $\Sp(A)$, we have:
    \begin{align*}
     \Perf(\Sp(A))\otimes_{\Perf(k)}\Perf(K) & \simeq \Perf(A)\otimes_{\Perf(k)}\Perf(K) \\
     &\simeq \Perf(A\otimes_{k}K) \\
     & \simeq \Perf(\Sp(A)\times_{\Sp(k)}\Sp(K)).   
    \end{align*}
    The first equivalence follows from \Cref{thm: Tate acyclicity}, and the last follows from the fact that, in this case, the tensor product $A\otimes_{k}K$ can be identified with the completed tensor product. Since $X$ is separated, we can take a cover by affinoid spaces and express $\Perf(X)$ as a limit of \cates of the form $\Perf(\Sp(A_{i}))$. We can then compute the tensor product as follows: 
    \begin{align*}
        \Perf(X) \otimes_{\Perf(k)}\Perf(K) & \simeq \lim \left( \Perf(\Sp(A_{i})) \otimes_{\Perf(k)}\Perf(K) \right) \\ 
        &\simeq \lim \Perf(\Sp(A_{i})\times_{\Sp(k)}\Sp(K)) \\ 
        &\simeq \Perf(X_K).
    \end{align*}
    The final statement follows from the fact that smooth and proper \cates are stable under base change.
\end{proof}

\begin{lemma}\label{lem: pb of the compact generator}
    Let $p\colon X \to \Sp(k)$ be a quasi-compact and separated rigid analytic variety such that $\mathrm{Ind}(\Perf(X))$ admits a compact generator $\m{E}$. Let $K$ be a finite field extension of $k$. We denote by $X_K$ the corresponding pullback rigid analytic variety: 
    \[
    \begin{tikzcd}
	{X_K} && X \\
	\\
	{\Sp(K)} && {\Sp(k).}
	\arrow[from=1-1, to=1-3, "i_{X}"]
	\arrow[from=1-1, to=3-1]
	\arrow[from=1-3, to=3-3]
	\arrow[from=3-1, to=3-3, "i"]
\end{tikzcd}
    \]
    Then $i_X^\ast(\m{E})$ is a compact generator of $\mathrm{Ind}(\Perf(X_K))$. Moreover, if $\Perf(X)$ is smooth and proper over $k$ and $B$ is the $\mathbb{E}_{1}$-algebra of endomorphisms of $\m{E}$, then $\Perf(X_K)$ is equivalent to $\LMod_{B\otimes_{k}K}(\Perf(K))$. In particular, the $K$-algebra of endomorphisms of $i^{\ast}_{X}(\m{E})$ is Morita equivalent to the $K$-algebra $B\otimes_{k}K$.
    \end{lemma}
\begin{proof}
    We observe that $i^{\ast}_{X}$ preserves perfect complexes, so $i_X^\ast(\m{E})$ is a compact object in $\mathrm{Ind}(\Perf(X_K))$. To show that $i^{\ast}_X(\m{E})$ is a generator, it suffices to prove that the right adjoint of $\mathrm{Ind}(i^{\ast}_{X})$ is conservative. We observe that $i_{X,\ast}(-)$ preserves perfect complexes. Indeed, being a perfect complex is a local property, allowing us to assume that $X$ is an affinoid space $\Sp(A)$. In this case, $X_K$ is identified with $\Sp(A\otimes_k K)$. Using \Cref{thm: Tate acyclicity}, restricted to perfect complexes, $i_{X,\ast}$ is identified with the forgetful functor, which preserves perfect complexes because $K$ is a finite extension of $k$. Thus, we can identify the right adjoint of $\mathrm{Ind}(i^{\ast}_{X})$ with $\mathrm{Ind}(i_{X,\ast})$. Furthermore, locally on $X$, we can identify $\mathrm{Ind}(i_{X,\ast})$ with the forgetful functor between module \cates, which is conservative. This implies that $i_{X}^{\ast}\m{E}$ is a generator. The final statement is deduced by first applying \Cref{lem: compact generation after bc} to conclude that $\Perf(X_K)$ is smooth and proper over $\Perf(K)$, and subsequently applying \Cref{lem: compact generator for perf and adjuntion}.
\end{proof}

\begin{lemma}\label{lem: diagramma back-chevalley}
     Consider the following diagram of algebras over $k$, where $B'=A'\otimes_{A}B$:
     \begin{equation*}
         \begin{tikzcd}
A \arrow[d, '] \arrow[r, ] & A' \arrow[d, ] \\
B \arrow[r, ']            & B'                    
\end{tikzcd}
     \end{equation*}
     For every $M\in \Mod_{A'}$, we can construct the following induced diagram:
\begin{equation*}
    \begin{tikzcd}
	{\Mod_{A}} &&& {\Mod_{A'}} \\
	\\
	\\
	{\Mod_{B}} &&& {\Mod_{B'}}
	\arrow[""{name=0, anchor=center, inner sep=0}, "{M\otimes_{A}-}"', shift right=3, from=1-1, to=1-4]
	\arrow["{B\otimes_{A}-}"', from=1-1, to=4-1]
	\arrow[""{name=1, anchor=center, inner sep=0}, shift right=2, from=1-4, to=1-1]
	\arrow["{B'\otimes_{A'}-}", from=1-4, to=4-4]
	\arrow[""{name=2, anchor=center, inner sep=0}, "{B'\otimes_{B}}"', shift right=3, from=4-1, to=4-4]
	\arrow[""{name=3, anchor=center, inner sep=0}, shift right=2, from=4-4, to=4-1]
	\arrow["\dashv"{anchor=center, rotate=90}, draw=none, from=0, to=1]
	\arrow["\dashv"{anchor=center, rotate=90}, draw=none, from=2, to=3]
\end{tikzcd}
\end{equation*}
 In this situation, the diagram formed by the left adjoints commutes. Furthermore, let $BC_{M}$ denote the Beck-Chevalley transformation associated with this diagram, and let $\m{C}$ be the full subcategory of $\Mod_{A'}$ consisting of modules $M$ for which $BC_{M}$ is an equivalence. Then $\Perf(A')\subseteq \m{C}$.
\end{lemma}
\begin{proof}
    We first note that $A'$ belongs to $\m{C}$. In this case, the right adjoint to the tensor product is simply the forgetful functor, and the Beck-Chevalley transformation is trivially an equivalence because for every module $N \in \Mod_{A'}$, the relation $B'=A'\otimes_{A}B$ implies:
 \[
 B\otimes_{A}N\simeq B'\otimes_{A'}N.
 \]
 To complete the proof, we note that if $M$ is a perfect complex over $A'$, it can be obtained as a finite colimit of retracts of finite direct sums of $A'$. Since the subcategory $\m{C}$ is closed under retracts and finite colimits, it follows that $M \in \m{C}$. 
\end{proof}

\begin{proof}[Proof of \Cref{prop: Perf vs AnM}]
    Fix a derived affinoid space $S$ with global sections algebra $A$. Our goal is to demonstrate a functorial equivalence of \cates between $\uline{\Perf}_{X}(S)$ and $\m{AnM}_{\Perf(X)}(S)$.
    Recall that since $\Perf(X)$ is smooth and proper, by \Cref{lem: compact generator for perf and adjuntion} it admits a compact generator $\m{E}$ and can be identified with $\Perf(B)$, where $B$ is the $\mathbb{E}_{1}$-algebra of endomorphisms of $\m{E}$. Using \Cref{rmk: AnM(B) and LModB}, we obtain an equivalence: 
    \[
    \m{AnM}_{\Perf(X)}(S) \coloneqq \Fun^{\mathrm{st}}_{k}(\Perf(X)^{\mathrm{op}}, \Perf(S)) \simeq \LMod_{B}(\Perf(S)).
    \]
    Furthermore, by definition:
    \[
    \uline{\Perf}_{X}(S) \coloneqq \Perf(X \times S).
    \]
    Consider the projection diagram:
    \begin{equation*}
    \begin{tikzcd}
	& {X \times S} \\
	{X} && {S}
	\arrow["{p_{1}}", from=1-2, to=2-1]
	\arrow["{p_{2}}"', from=1-2, to=2-3]
    \end{tikzcd}
    \end{equation*}
    Using the compact generator, we can establish an adjunction: 
    \begin{equation}\label{eq: 1}
    \begin{tikzcd}
	{\Perf(X\times S)} &&&& {\LMod_{B}(\Perf(S))}.
	\arrow[""{name=0, anchor=center, inner sep=0}, "{\psi\coloneqq p_{2,\ast}\m{Hom}_{X\times S}(p_1^{\ast}\m{E},-)}"', shift right=2, from=1-1, to=1-5]
	\arrow[""{name=1, anchor=center, inner sep=0}, "{\phi\coloneqq p_2^{\ast}(-) \otimes p_1^{\ast}\m{E}}"', shift right=2, from=1-5, to=1-1]
	\arrow["\dashv"{anchor=center, rotate=-90}, draw=none, from=1, to=0]
    \end{tikzcd}
    \end{equation}
    By \Cref{lem: compact generator for perf and adjuntion}, the adjunction in \eqref{eq: 1} is an equivalence when $S\simeq \Sp(k)$, as these are exactly the functors producing the adjoint equivalence between $\Perf(X)$ and $\Perf(B)$. Using \Cref{lem: compact generator for perf and adjuntion}, \Cref{lem: compact generation after bc}, and \Cref{lem: pb of the compact generator}, we deduce that this adjunction is also an equivalence when $S \simeq \Sp(K)$ for every finite field extension $K$ of $k$.\\
    The strategy is to deduce the general case from the fields $S\simeq \Sp(K)$. For every closed point $s \colon \Sp(K) \to S$ with residue field $K$, we consider the Cartesian diagram: 
    \begin{equation}\label{eq: diagram notazioni 1}
\begin{tikzcd}
	{X_K\coloneq X\times \Sp(K)} && {X \times S} \\
	\\
	{\Sp(K)} && S
	\arrow["{i_{s}}", from=1-1, to=1-3]
	\arrow[from=1-1, to=3-1, "p_{s}"]
	\arrow[from=1-3, to=3-3, "p_2"]
	\arrow["s"', from=3-1, to=3-3]
\end{tikzcd}
    \end{equation}
    This induces the following diagram: 
    
    \begin{equation}\label{eq: diagram right adjointable}
        \begin{tikzcd}
	{\Perf(X\times S)} &&& {\LMod_{B}(\Perf(S))} & \\
	\\
	{\Perf(X_K)} &&& {\LMod_{B}(\Perf(K)) \mathrlap{\simeq \Perf(B\otimes_{k}K)}} & {}
	\arrow[""{name=0, anchor=center, inner sep=0}, "\phi"', shift right=3, from=1-1, to=1-4]
	\arrow["{i_{s}^{\ast}}"', from=1-1, to=3-1]
	\arrow[""{name=1, anchor=center, inner sep=0}, "\psi"', shift right=2, from=1-4, to=1-1]
	\arrow["{s^{\ast}}", from=1-4, to=3-4]
	\arrow[""{name=2, anchor=center, inner sep=0}, shift right=3, from=3-1, to=3-4]
	\arrow[""{name=3, anchor=center, inner sep=0}, shift right=2, from=3-4, to=3-1]
	\arrow["\dashv"{anchor=center, rotate=-90}, draw=none, from=1, to=0]
	\arrow["\dashv"{anchor=center, rotate=-90}, draw=none, from=3, to=2]
\end{tikzcd}
    \end{equation}
We observe that this diagram commutes when traversing the left adjoints. 
 Because derived affinoids are Jacobson, the family of pullback functors $\{i_s^{\ast}\}_{s\in S}$ is jointly conservative. Furthermore, we know that the bottom adjunction is an equivalence. Therefore, to prove that the unit $\eta$ and counit $\epsilon$ of the adjunction between $\phi$ and $\psi$ are equivalences, it suffices to show that the map $i_s^{\ast}$ maps the unit of the top adjunction to the unit of the bottom adjunction, and similarly for the counit via $s^{\ast}$. This reduces to verifying that the diagram is right adjointable \cite[Definition 4.7.4.13]{HA}. 
 
 By \Cref{lem: compact generation after bc}, the object $i^{\ast}_s p^{\ast}_{1}(\m{E})$ is a generator of $\Perf(X_K)$. Since the forgetful functor from $\LMod_{B}(\Perf(K))$ to $\Perf(K)$ is conservative, the right adjointability of the diagram in \eqref{eq: diagram right adjointable} is equivalent to the right adjointability of the following expanded diagram:
\begin{equation}
\begin{tikzcd}
	{\Perf(X\times S)} &&& {\Perf(X\times S)} &&& {\Perf(S)} \\
	\\
	{\Perf(X_K)} &&& {\Perf(X_K)} & {} && {\Perf(K)}
	\arrow[""{name=0, anchor=center, inner sep=0}, "{\m{Hom}_{X\times S}(p_{1}^{\ast}\m{E},-)}"', shift right=3, from=1-1, to=1-4]
	\arrow["{i_{s}^{\ast}}"', from=1-1, to=3-1]
	\arrow[""{name=1, anchor=center, inner sep=0}, "{-\otimes p^{\ast}_{1}\m{E}}"', shift right=2, from=1-4, to=1-1]
	\arrow[""{name=2, anchor=center, inner sep=0}, "{p_{2,\ast}}"', shift right=3, from=1-4, to=1-7]
	\arrow["{i_{s}^{\ast}}", from=1-4, to=3-4]
	\arrow[""{name=3, anchor=center, inner sep=0}, "{p^{\ast}_{2}}"', shift right=2, from=1-7, to=1-4]
	\arrow["{s^{\ast}}", from=1-7, to=3-7]
	\arrow[""{name=4, anchor=center, inner sep=0}, "{\m{Hom}_{X_K}(i^{\ast}_{s}p^\ast_1\m{E},-)}"', shift right=3, from=3-1, to=3-4]
	\arrow[""{name=5, anchor=center, inner sep=0}, "{-\otimes i^{\ast}_{s}p^{\ast}_1\m{E}}"', shift right=2, from=3-4, to=3-1]
	\arrow[""{name=6, anchor=center, inner sep=0}, "{p_{s,\ast}}"', shift right=3, from=3-4, to=3-7]
	\arrow[""{name=7, anchor=center, inner sep=0}, "{p^{\ast}_{s}}"', shift right=2, from=3-7, to=3-4]
	\arrow["\dashv"{anchor=center, rotate=-90}, draw=none, from=1, to=0]
	\arrow["\dashv"{anchor=center, rotate=-90}, draw=none, from=3, to=2]
	\arrow["\dashv"{anchor=center, rotate=-90}, draw=none, from=5, to=4]
	\arrow["\dashv"{anchor=center, rotate=-90}, draw=none, from=7, to=6]
\end{tikzcd}
\end{equation}
We prove this by showing that both the left and right squares are individually right adjointable. The right adjointability of the square on the right follows from the base change property applied to the diagram in \eqref{eq: diagram notazioni 1}, which was proven in \cite[Lemma 6.4]{derhom}. The right adjointability of the square on the left follows utilizing the fact that perfect complexes are dualizable and the functor $i^{\ast}_{s}$ is symmetric monoidal.

To prove that this equivalence restricts to simple objects, for every $\m{F} \in \Perf(X\times S)$ we consider the following diagram:
\begin{equation}
\begin{tikzcd}
	&& {\Perf(S)} && \\
	&& {} \\
	\\
	{\Perf(X\times S)} &&&& {\LMod_{B}(\Perf(S)).}
	\arrow[""{name=0, anchor=center, inner sep=0}, "{\m{F}\otimes-}", shift left, from=1-3, to=4-1]
	\arrow[""{name=1, anchor=center, inner sep=0}, "{\psi(\m{F})\otimes-}"', shift right, from=1-3, to=4-5]
	\arrow[""{name=2, anchor=center, inner sep=0}, "{\uline{\m{Hom}}_{X\times S}(\m{F},-)}", shift left=3, from=4-1, to=1-3]
	\arrow[""{name=3, anchor=center, inner sep=0}, "{\psi\coloneqq p_{1\ast}\m{Hom}_{X\times S}(p^{\ast}_{1}\m{E},-)}"', shift right=2, from=4-1, to=4-5]
	\arrow[""{name=4, anchor=center, inner sep=0}, "{\uline{\Hom}_{B}(\psi(\m{F}),-)}"', shift right=3, from=4-5, to=1-3]
	\arrow[""{name=5, anchor=center, inner sep=0}, "{\phi\coloneqq p^{\ast}_{1}\m{E}\otimes p_{2}^{\ast}(-)}"', shift right=2, from=4-5, to=4-1]
	\arrow["\dashv"{anchor=center, rotate=45}, draw=none, from=1, to=4]
	\arrow["\dashv"{anchor=center, rotate=138}, draw=none, from=0, to=2]
	\arrow["\dashv"{anchor=center, rotate=-90}, draw=none, from=5, to=3]
\end{tikzcd}
\end{equation}
We adopt the following notation: 
\[
L\coloneqq \psi(\m{F})\otimes- \; ; \; L'\coloneqq \m{F}\otimes- \; ; \; R\coloneqq \uline{\Hom}_{B}(\psi(\m{F}),-) \; ; \; R'\coloneqq \uline{\m{Hom}}_{X\times S}(\m{F},-).
\]
The diagram commutes on the left adjoints. Indeed, we have:
\begin{align*}
    \phi\circ L & \coloneqq \phi(\psi(\m{F})\otimes-)\\ 
    & \simeq p^{\ast}_{1}\m{E}\otimes p^{\ast}_{2}(\psi(\m{F})\otimes -) \\ 
    & \simeq p^{\ast}_{1}(\m{E})\otimes p^{\ast}_{2}(\psi\m{F}) \otimes p^{\ast}_{2}(-) \\ 
    & \simeq \phi(\psi\m{F})\otimes p^{\ast}_{2}(-) \\ 
    & \simeq \m{F}\otimes p^{\ast}_{2}(-) \coloneqq L'.
\end{align*}
Since the adjunction between $\phi$ and $\psi$ is an equivalence, we deduce that 
\[
L \simeq \psi \circ L'.
\]
Passing to the right adjoints implies that 
\[
R \simeq R' \circ \phi.
\]
Evaluating at the respective identity objects, we obtain: 
\[
\uline{\Hom}_{B}(\psi(\m{F}), \psi(\m{F})) \simeq \uline{\m{Hom}}_{X\times S}(\m{F},\m{F}).
\]
Consequently, $\m{F}$ is simple in the sense of \Cref{def: simple analytic} if and only if $\psi(\m{F})$ is simple in the sense of \Cref{def: simple affine setting}. This proves that the equivalence $\psi$ restricts to an equivalence between simple objects.
\end{proof}

In the following, we denote by $t_{0}\uline{\Perf}^{\mathrm{sim}}_{X}$ the restriction of $\uline{\Perf}^{\mathrm{sim}}_{X}$ to underived affinoids. Note that by \Cref{cor: simple analytic underived}, it is naturally a rigid analytic $1$-stack. 

\begin{corollary}\label{cor: analytifiction of the trouncation}
    Let $X \to \Sp(k)$ be a proper, quasi-compact, and separated rigid analytic variety such that $\Perf(X)$ is a smooth and proper \cat over $\Perf(k)$. Then the analytic $1$-stacks 
    $t_{0}\uline{\Perf}^{\mathrm{sim}}_{X}$ and $(t_{0}\m{M}^{\mathrm{sim}}_{\Perf(X)})^{\mathrm{an}}$ are equivalent.
\end{corollary}
\begin{proof}
Combining \Cref{prop: MB sim vs AnMB sim} and \Cref{prop: Perf vs AnM}, we obtain an equivalence between $\uline{\Perf}^{\mathrm{sim}}_{X}$ and $(\m{M}^{\mathrm{sim}}_{\Perf(X)})^{\mathrm{an}}$. Moreover, since the analytification functor preserves geometric stacks (see \cite[Lemma 2.24]{porta2016higher}), we deduce that $(t_{0}\m{M}^{\mathrm{sim}}_{\Perf(X)})^{\mathrm{an}}$ is a rigid analytic geometric $1$-stack. Using \Cref{thm: teorema mauro} and \Cref{thm: Tate acyclicity}, the equivalence in \Cref{prop: MB sim vs AnMB sim} restricts to an equivalence between the rigid analytic $1$-stacks $(t_0\m{M}^{\mathrm{sim}}_{\Perf(X)})^{\mathrm{an}}$ and $t_0\m{AnM}^{\mathrm{sim}}_{\Perf(X)}$. We conclude by applying \Cref{prop: Perf vs AnM}, observing that all the stacks considered are naturally analytic $1$-stacks, as explained above.
\end{proof}

\section{Main theorem}
The aim of this section is to prove the following main result:

\begin{theorem}\label{thm: main}
    Let $X \to \Sp(k)$ be a smooth and proper, quasi-compact, and separated rigid analytic variety. Then $\Perf(X)$ is smooth and proper over $\Perf(k)$ if and only if $X$ is algebraizable (i.e., there exists an algebraic space $Y$ and an equivalence $Y^{\mathrm{an}} \simeq X$).
\end{theorem}

In the final subsection, we shall also deduce an analogous statement for formal schemes. 

One direction of the proof is straightforward, as an algebraic space $Y$ is smooth and proper if and only if $\Perf(Y)$ is smooth and proper (\cite[Theorem 11.4.0.3]{SAG}). Furthermore, by GAGA theorems (see \cite{porta2016higher}), we obtain the same result for $\Perf(Y^{\mathrm{an}})$. In the formal setting, an analogous argument applies using the Grothendieck Existence Theorem (\cite[Theorem 8.5.0.3]{SAG}).

The converse direction is more involved. We provide here an overview of the proof strategy, with the technical details deferred to the subsequent subsections. We begin by observing that, by \Cref{lem: compact generator for perf and adjuntion}, there exists an equivalence between $\Perf(X)$ and $\Perf(B)$, where $B$ is the $k$-algebra of endomorphisms of a compact generator of $\Perf(X)$. Moreover, by \Cref{cor: categories vs algebras}, $B$ is a smooth and proper $k$-algebra. 

We consider the truncation of the derived algebraic stack $\m{M}^{\mathrm{sim}}_{\Perf(X)}$, denoted by $t_{0}\m{M}^{\mathrm{sim}}_{\Perf(X)}$. According to \cite{modg}, this is an Artin $1$-stack that admits a good moduli space $\mathrm{M}_{X}$. The projection map
\begin{equation}\label{eq: coarse moduli}
t_{0}\m{M}^{\mathrm{sim}}_{\Perf(X)} \to \mathrm{M}_{X}    
\end{equation}
is a $\mathbb{G}_{m}$-gerbe. In particular, $\mathrm{M}_{X}$ is an algebraic space locally of finite type, and its analytification $\mathrm{M}_{X}^{\mathrm{an}}$ is an underived rigid analytic space.

The proof then proceeds via the following main steps: 
\begin{enumerate}
    \item We prove that $X$ can be embedded as an open subspace into $\mathrm{M}^{\mathrm{an}}_{X}$.
    \item We establish an algebraization result for such an embedding. 
\end{enumerate}

The proofs of these two points will be detailed in the next subsections, completing the proof of the main theorem. Before proceeding, we verify that the properness of $\Perf(X)$ holds whenever $X$ is proper; however, as \Cref{thm: main} suggests, this correspondence does not extend to smoothness in the same way.

\begin{corollary}\label{cor: non commutative properness}
    Let $X$ be a proper, quasi-compact rigid analytic variety. Then $\Perf(X)$ is a proper $k$-linear \cat.
\end{corollary}
\begin{proof}
    We must show that $\Hom_{X}(\m{F},\m{G})$ is a perfect $k$-module for any pair of perfect complexes $\m{F}$ and $\m{G}$. This follows directly from \Cref{lem: non commutative properness} and \Cref{lemma: hom perfect}.
\end{proof}

\begin{remark}
	Note that by the previous corollary, to prove algebraizability it suffices to verify the smoothness of $\Perf(X)$. This property represents the key categorical distinction between the algebraic and analytic settings.
\end{remark}

\subsection{Embedding into the analytification of the algebraic space}

Let $X \to \Sp(k)$ be a smooth and proper, quasi-compact, and separated rigid analytic variety, such that $\Perf(X)$ is smooth and proper over $\Perf(k)$. We denote by $\mathrm{M}_X$ the coarse moduli space of $t_0\m{M}^{\mathrm{sim}}_{\Perf(X)}$, as established in \eqref{eq: coarse moduli}.

In this subsection, we construct a map 
\begin{equation}\label{eq: mappa verso l'analytificazione del moduli space}
    \alpha \colon X \times \mathrm{B}\mathbb{G}^{\mathrm{an}}_m \to t_0\uline{\Perf}^{\mathrm{sim}}_{X}\simeq (t_0\m{M}_{\Perf(X)})^{\mathrm{an}}.
\end{equation}
Subsequently, we will prove the following proposition:

\begin{proposition}\label{prop: X open into the analytification}
    Under the assumptions of \Cref{thm: main}, the map constructed in \eqref{eq: mappa verso l'analytificazione del moduli space} induces an open immersion of $X$ into $\mathrm{M}_X^{\mathrm{an}}$.
\end{proposition}

To construct the map in \Cref{eq: mappa verso l'analytificazione del moduli space}, we begin with the following setup. The functor $\mathcal{X}$ introduced below can be viewed as the rigid analytic analogue of the algebraic stack of coherent sheaves of length $1$, as studied in \cite{anel2007derivedequivalenceclassesalgebraic}.

\begin{construction}
    Consider the following (underived) functor:
    \[
    \mathcal{X}\colon \Afd_{k} \to \mathcal{S}
    \]
    which assigns to an affinoid space $S$ the core of the full subcategory of $\Perf(X\times S)$ spanned by objects $\m{F}$ having Tor-amplitude $[0,0]$ relative to $S$ (i.e., those which are flat over $S$, and whose only non-zero cohomology sheaf is coherent and concentrated in degree $0$). Furthermore, for every closed point $s \in S$, we require the fiber $\m{F}_{s}$ to be the skyscraper sheaf of length $1$ supported on a single point of the fiber $X_s$. 
    
    Observe that since the objects in $\mathcal{X}(S)$ lie in the heart of the t-structure on $\Mod_{\mathscr{O}_{X\times S}}$, their mapping spaces are $0$-truncated \cite[Remark 1.2.1.12]{HA}. In particular, $\mathcal{X}$ can be viewed as a functor taking values in $1$-groupoids.
\end{construction}

\begin{lemma}\label{lem: equivalenza con X x BGm}
    There is an equivalence of presheaves between $\mathcal{X}$ and $X \times \mathrm{B}\mathbb{G}^{\mathrm{an}}_{m}.$
\end{lemma}

\begin{proof}
    Recall that the functor $\mathcal{X}$ takes values in $1$-groupoids, and the same holds for the product stack $X \times \mathrm{B}\mathbb{G}^{\mathrm{an}}_m$. It therefore suffices to fix an arbitrary (classical) affinoid space $S$ and construct an equivalence of $1$-groupoids evaluated at $S$. 
    
    We begin by constructing the functor 
    \[
    \Psi \colon (X \times \mathrm{B}\mathbb{G}^{\mathrm{an}}_m)(S) \to \mathcal{X}(S).
    \]
    An object in the domain is a pair $(f, \mathcal{L})$, where $f \colon S \to X$ is a morphism of analytic spaces and $\mathcal{L}$ is a line bundle on $S$. Considering the graph morphism $\Gamma_f \colon S \to X \times S$, we define:
    \[
        \Psi(f, \mathcal{L}) \coloneqq \Gamma_{f_{*}} \mathcal{L}.
    \]
    Observe that, since $X$ is separated, $\Gamma_f$ is a closed immersion and $\mathcal{L}$ is locally free over $S$. We verify the relative Tor-amplitude condition via the projection formula \cite[Theorem 5.7]{derhom}. For any discrete $\mathscr{O}_S$-module $N$, we compute the derived tensor product:
    \[
        (\Gamma_f)_* \mathcal{L} \otimes_{\mathscr{O}_{X \times S}} p_S^* N \simeq \Gamma_{f_*} \left( \mathcal{L} \otimes_{\mathscr{O}_S} \Gamma_f^* p_S^* N \right).
    \]
    Because $p_S \circ \Gamma_f \simeq \mathrm{id}_S$, we have $\Gamma_f^* p_S^* N \simeq N$, and the expression simplifies to $\Gamma_{f_*} (\mathcal{L} \otimes_{\mathscr{O}_S} N)$. Since $\mathcal{L}$ is a line bundle on $S$ (in particular, it is locally free and concentrated in degree $0$), the derived tensor product $\mathcal{L} \otimes^{\mathbb{L}}_{\mathscr{O}_S} N$ is also concentrated in degree $0$. Finally, since $\Gamma_f$ is a closed immersion, its pushforward is t-exact by \cite[Lemma 6.2]{derhom}. Thus, the resulting complex $\Gamma_{f_*} \mathcal{L}$ has Tor-amplitude $[0,0]$ relative to $S$.
    
    Furthermore, we must verify the condition on the fibers. Let $s \in S$ be a closed point, and consider the inclusion of the fiber $j_s \colon X_s \to X \times S$. To compute $j_s^* \Gamma_{f_*} \mathcal{L}$, we consider the Cartesian square:
    \begin{equation}\label{eq: base change comparison fucntor}
    \begin{tikzcd}
        \mathrm{Sp}(k(s)) \arrow[r, "i_s"] \arrow[d, "i_{\tilde{x}}"'] & S \arrow[d, "\Gamma_f"] \\
        X_s \arrow[r, "j_s"'] & X \times S
    \end{tikzcd}
    \end{equation}
    where $i_s$ is the standard inclusion of the point $s$, and the vertical morphism $i_{\tilde{x}}$ is the inclusion of a $k(s)$-rational point $\tilde{x} \in X_s$. We observe that $\tilde{x}$ is a lift of the point $f(s) \in X$ to $X_s$.
    
    Since $\Gamma_f$ is a proper morphism (being a closed immersion), the derived base change theorem \cite[Theorem 6.8]{derhom} yields a canonical isomorphism:
    \[
        j_s^* \Gamma_{f_*} \mathcal{L} \simeq i_{\tilde{x}_*} i_s^* \mathcal{L}.
    \]
    Since $\mathcal{L}$ is a line bundle, its pullback $i_s^* \mathcal{L}$ to the point $s$ is a $1$-dimensional vector space over $k(s)$ concentrated in degree $0$. Moreover, since $i_{\tilde{x}}$ is a closed immersion, $i_{\tilde{x}_*} i_s^* \mathcal{L}$ is a discrete skyscraper sheaf on $X_s$ supported at the point $\tilde{x}$ of length $1$, perfectly matching the functor's requirements.
    
    Regarding morphisms, since $\mathcal{X}(S)$ is a $1$-groupoid, a morphism
    \[
    \alpha \colon (f,\mathcal{L}) \to (f',\mathcal{L}')
    \]
    is given by a map $f \to f'$ and an isomorphism of line bundles $\mathcal{L} \xrightarrow{\sim} \mathcal{L}'$. Because the mapping space $\mathrm{Map}(S,X)$ is $0$-truncated, any such map implies $f \simeq f'$. Thus, it suffices to consider the case where $f = f'$, allowing us to define 
    \[
    \Psi(\alpha) \coloneqq \Gamma_{f_*}(\alpha).
    \]

    We now prove that $\Psi$ is fully faithful and essentially surjective. We begin by demonstrating fully faithfulness. 
    Let $(f_1, \mathcal{L}_1)$ and $(f_2, \mathcal{L}_2)$ be two objects in $(X \times \mathrm{B}\mathbb{G}^{\mathrm{an}}_m)(S)$. We must show that $\Psi$ induces a bijection on the associated mapping spaces. As explained above, the mapping space in the domain is empty if $f_1 \neq f_2$, and it naturally identifies with $\mathrm{Isom}_S(\mathcal{L}_1, \mathcal{L}_2)$ if $f_1 = f_2 = f$. 

    Suppose first that $f_1 \neq f_2$. The sheaves $\Gamma_{f_{1_*}} \mathcal{L}_1$ and $\Gamma_{f_{2_*}} \mathcal{L}_2$ are topologically supported on the closed subsets $\Gamma_{f_1}(S)$ and $\Gamma_{f_2}(S)$ in $X \times S$. Since $X$ is separated, the graph of any morphism is a closed immersion, and for distinct morphisms $f_1 \neq f_2$, their graphs are distinct subspaces. This implies 
    \[
    \mathrm{Isom}_{\mathcal{X}(S)}(\Gamma_{f_{1*}} \mathcal{L}_1, \Gamma_{f_{2*}} \mathcal{L}_2) = \emptyset,
    \] 
    which exactly matches the empty mapping space of the domain.

    Suppose now that $f_1 = f_2 = f$. We must verify that the pushforward functor along the graph induces a bijection:
    \[
    \Gamma_{f_*} \colon \mathrm{Isom}_S(\mathcal{L}_1, \mathcal{L}_2) \xrightarrow{\sim} \mathrm{Isom}_{\mathcal{X}(S)}(\Gamma_{f_*} \mathcal{L}_1, \Gamma_{f_*} \mathcal{L}_2).
    \]
    Since $p_S \circ \Gamma_f = \mathrm{id}_S$, injectivity is immediate: if $\Gamma_{f_*}(\gamma_1) = \Gamma_{f_*}(\gamma_2)$, applying $p_{S_*}$ to both sides yields $p_{S_*}\Gamma_{f_*}(\gamma_1) = p_{S_*}\Gamma_{f_*}(\gamma_2)$, which simplifies directly to $\gamma_1 = \gamma_2$. 

    To demonstrate surjectivity, let $\beta \colon \Gamma_{f_*} \mathcal{L}_1 \xrightarrow{\sim} \Gamma_{f_*} \mathcal{L}_2$ be an arbitrary isomorphism in $\mathcal{X}(S)$. We define a candidate isomorphism in the domain as $\gamma \coloneqq p_{S_*}(\beta) \colon \mathcal{L}_1 \xrightarrow{\sim} \mathcal{L}_2$. We must show that $\Gamma_{f_*}(\gamma) \simeq \beta$.
    Observe that both $\beta$ and $\Gamma_{f_*}(\gamma)$ are supported on the graph of $f$. Furthermore, the restriction of $\Gamma_f \circ p_{S}$ to the graph of $f$ is equivalent to the identity functor. This implies that 
    \[
    \Gamma_{f_*} \circ p_{S_*} (\beta) \simeq \beta,
    \]
    concluding the proof that $\Psi$ is fully faithful. 

    We now turn to the essential surjectivity of $\Psi$. Let $\mathscr{F}$ be an object in $\mathcal{X}(S)$; we will construct a pair $(f,\mathcal{L})$ in $(X \times \mathrm{B}\mathbb{G}^{\mathrm{an}}_m)(S)$ such that $\Psi(f,\mathcal{L}) \simeq \mathscr{F}$.

    By definition, the sheaf $\mathscr{F}$ in $\mathcal{X}(S)$ corresponds to an object in $\Perf(X\times S)$ such that for every closed point $s$ in $S$, the fiber $\mathscr{F}_s$ is supported at a single point in $X_s$. Consequently, the support of $\mathscr{F}$ can be described as the graph of the function 
    \[
    \mathrm{supp}(\mathscr{F})(-)\colon S \to X
    \]
    sending a point $s$ to the projection onto $X$ of the single point in the fiber $X_s$ that supports $\mathscr{F}_s$. Moreover, since $\mathscr{F}$ has Tor-amplitude $[0,0]$ relative to $S$ and $X$ is proper, we deduce that
    \[
        \mathcal{L} \coloneqq p_{S_*} \mathscr{F}
    \]
    is a vector bundle on $S$. The specific requirement on the fiber $\mathscr{F}_s$ ensures that $\mathcal{L}$ is, in fact, a line bundle on $S$.

    We need to show that 
    \[
    \Psi(\mathrm{supp}(\mathscr{F})(-), \mathcal{L}) \simeq \mathscr{F}.
    \]
    We observe that the support of $\mathscr{F}$ coincides with the graph of the function $\mathrm{supp}(\mathscr{F})(-)$. In particular, the composition $\Gamma_{\mathrm{supp}(\mathscr{F})(-)} \circ p_S$, when restricted to this graph, is equivalent to the identity functor. Therefore, we obtain 
    \[
    \Psi(\mathrm{supp}(\mathscr{F})(-), \mathcal{L}) \coloneqq \Gamma_{\mathrm{supp}(\mathscr{F})(-)_{*}} (p_{S_*}(\mathscr{F})) \simeq \mathscr{F}.
    \]
    This completes the proof.
\end{proof}

As a consequence of \Cref{lem: equivalenza con X x BGm}, $\mathcal{X}$ is a stack. Furthermore, we observe that the objects of $\mathcal{X}(S)$ are simple sheaves over $S$ in the sense of \Cref{def: simple analytic}. In particular, this induces a canonical map:
\begin{equation}\label{eq: mappa inclusione}
    \tilde{\alpha} \colon \mathcal{X} \to t_0\uline{\Perf}_{X}^{\mathrm{sim}}.
\end{equation}

\begin{lemma}\label{lem: open immersion}
    The map constructed in \eqref{eq: mappa inclusione} is an open immersion.
\end{lemma}
\begin{proof}
    We will prove that the canonical map \eqref{eq: mappa inclusione} is a monomorphism and formally \'etale. Note that, by \Cref{cor: analytifiction of the trouncation}, these stacks are locally of finite presentation, so this will imply it is an open immersion. 
    
    By definition, for every affinoid space $S$, the \cat $\mathcal{X}(S)$ is a full subcategory of $\Perf(X\times S)^{\mathrm{sim}}$. This immediately implies that the diagonal map
    \[
    \mathcal{X} \to \mathcal{X}\times_{t_{0}\uline{\Perf}^{\mathrm{sim}}_{X}} \mathcal{X}
    \]
    is an equivalence, hence $\tilde{\alpha}$ is a monomorphism. 
    
    We now prove that $\tilde{\alpha}$ is formally \'etale. Let $\Sp(A)$ be an affinoid space, $I$ a square-zero ideal in $A$, and consider the following lifting problem:
\[\begin{tikzcd}
	{\mathcal{X}} && {\Sp(A/I)} \\
	\\
	{t_0\uline{\Perf}^{\mathrm{sim}}_{X}} && {\Sp(A).}
	\arrow["{\tilde{\alpha}}"', from=1-1, to=3-1]
	\arrow[from=1-3, to=1-1]
	\arrow[from=1-3, to=3-3]
	\arrow[dashed, from=3-3, to=1-1]
	\arrow[from=3-3, to=3-1]
\end{tikzcd}\]
    We must prove the existence of the dotted arrow. The above diagram is specified by a perfect complex $\m{E}$ on $X \times \Sp(A)$ such that, under the canonical inclusion
    \[
    j\colon X\times \Sp(A/I) \to X \times \Sp(A),
    \]
    the pullback $j^\ast \m{E}$ belongs to the subcategory $\mathcal{X}(\Sp(A/I))$. We need to show that $\m{E}$ itself lies in $\mathcal{X}(\Sp(A))$. 
    
    To prove that $\m{E}$ has Tor-amplitude $[0,0]$ relative to $\Sp(A)$, we can work locally and assume $X = \Sp(R)$ for some affinoid algebra $R$. Using \Cref{thm: Tate acyclicity} and \cite[Corollary 2.7.3.2]{SAG}, we deduce that $\m{E}$ indeed has Tor-amplitude $[0,0]$ relative to $\Sp(A)$. Finally, the condition that the restriction of $\m{E}$ to the fiber over every closed point in $\Sp(A)$ is a skyscraper sheaf is automatically satisfied: since $I$ is a nilpotent ideal, the affinoid spaces $\Sp(A)$ and $\Sp(A/I)$ share the exact same underlying topological space and closed points.  
\end{proof}

According to \cite[Corollary 3.22]{modg}, there exists an algebraic space locally of finite type, $\mathrm{M}_X$, such that the Artin stack $t_0\m{M}^{\mathrm{sim}}_{\Perf(X)}$ admits it as a good moduli space:
\begin{equation}\label{eq: mappa G_m gerbe}
p\colon t_0\m{M}^{\mathrm{sim}}_{\Perf(X)} \to \mathrm{M}_{X}.    
\end{equation}
Moreover, the map $p$ is a $\mathbb{G}_{m}$-gerbe. By analyzing the analytification of this map, we can conclude the proof of \Cref{prop: X open into the analytification}. 

\begin{proof}[Proof of \Cref{prop: X open into the analytification}]
    Taking the analytification of the map in \eqref{eq: mappa G_m gerbe} and invoking the equivalence established in \Cref{cor: analytifiction of the trouncation}, we obtain the following $\mathbb{G}^{\mathrm{an}}_{m}$-gerbe projection: 
    \[
    p^{\mathrm{an}}\colon t_0\uline{\Perf}^{\mathrm{sim}}_{X} \to \mathrm{M}^{\mathrm{an}}_{X}.
    \]
    Using \Cref{lem: equivalenza con X x BGm} and \Cref{lem: open immersion}, we can identify $X \times \mathrm{B}\mathbb{G}^{\mathrm{an}}_{m}$ as an open substack of $t_0\uline{\Perf}^{\mathrm{sim}}_{X}$. Since $p^{\mathrm{an}}$ is a $\mathbb{G}^{\mathrm{an}}_{m}$-gerbe, taking the quotient of this open substack by the gerbe structure exactly yields $X$. Because the projection map of a gerbe is smooth (and thus open), the image of this substack under $p^{\mathrm{an}}$ is canonically isomorphic to $X$, and it embeds as an open analytic subspace of $\mathrm{M}^{\mathrm{an}}_{X}$. 
\end{proof}

\subsection{Non-Archimedean Moishezon spaces}
The notion of non-Archimedean Moishezon spaces was developed and studied in \cite{conrad}. In this section, we recall the definition of non-Archimedean Moishezon spaces and state the main theorem characterizing them, from which we will deduce the algebraization result needed to complete the proof in the non-Archimedean case.\\

\begin{definition}[\cite{conrad}]\label{def: Moishezon non-arch}
    A rigid analytic variety $X$ is \emph{Moishezon} if it is proper and each irreducible component $X_{i}$ of $X$, endowed with its reduced structure, satisfies $\mathrm{trdeg}_{k}(K(X_{i})) = \dim X_{i}$, where $K(X_{i})$ is the field of meromorphic functions and $\mathrm{trdeg}_{k}$ denotes its transcendence degree over $k$.
\end{definition}

\begin{theorem}[\cite{conrad}]\label{thm: caratterizzazione Moishezon}
    The full subcategory of Moishezon spaces inside the category of rigid analytic varieties coincides with the essential image of the analytification functor from the category of proper algebraic spaces over $k$. 
\end{theorem}

The proof of the following corollary is inspired by \cite[Proposition 5.7]{tova}.

\begin{corollary}\label{cor: algebraization of opens}
    Let $M$ be a quasi-separated algebraic space, and let $X$ be a rigid analytic variety satisfying the assumptions of \Cref{thm: main}. If $j\colon X \to M^{\mathrm{an}}$ is an open immersion, then $X$ is algebraizable. 
\end{corollary}
\begin{proof}
   We will prove that $X$ is a Moishezon space.
   
   Since $X$ is quasi-compact, its image $j(X) \subset M^{\mathrm{an}}$ is also quasi-compact. Because $M$ is a quasi-separated algebraic space, it admits an \'etale covering by affine schemes. We can select a finite number of these affine schemes whose images cover $j(X)$. The union of these images forms a quasi-compact, open algebraic subspace $M' \subset M$ of finite type over $k$. The map $j$ then factors through an open immersion $X \to (M')^{\mathrm{an}}$. Consequently, we may assume without loss of generality that $M$ is quasi-compact and locally of finite type over $k$.
   
   To check the Moishezon condition, we evaluate $X$ component by component; thus, we may assume that $X$ is irreducible. The open immersion $j$ then factors through an open immersion $i\colon X \to Z$, where $Z$ is an irreducible component of $M^{\mathrm{an}}$. Furthermore, there exists an irreducible component $M_{0}\subset M$ such that $Z=M_{0}^{\mathrm{an}}$. Therefore, we may also assume that $M$ is irreducible. 
   
   Since $X$ and $M^{\mathrm{an}}$ are irreducible, their rings of meromorphic functions $K(M^{\mathrm{an}})$ and $K(X)$ are fields. The open immersion $j$ induces an injective field extension:
   \[
   K(M^{\mathrm{an}}) \hookrightarrow K(X).
   \]
   Let $n$ be the dimension of $X$. Because $j$ is an open immersion, we have
   \[
   \dim M^{\mathrm{an}} = \dim M = n.
   \]
   Since $M$ is an irreducible algebraic space of finite type, it contains a dense open affine subscheme $U$ of dimension $n$. We therefore have the field inclusions:
   \[
    K(U) \simeq K(M) \subset K(M^{\mathrm{an}}) \subset K(X).
   \]
   Since $U$ is an affine variety of dimension $n$ over $k$, its function field satisfies $\mathrm{trdeg}_{k}(K(U)) = n$. This implies that
   \[
   \mathrm{trdeg}_{k}(K(X)) \geq n.
   \]
   On the other hand, by \cite[Theorem 2.2.1]{conrad}, for any proper rigid analytic variety, we naturally have $\mathrm{trdeg}_{k}(K(X)) \leq \dim X = n$. We conclude that $\mathrm{trdeg}_{k}(K(X)) = n$, which means $X$ is a Moishezon space. The algebraizability of $X$ now follows directly from \Cref{thm: caratterizzazione Moishezon}.
\end{proof}

For completeness, we now assemble the complete proof of the main theorem.

\begin{proof}[Proof of \Cref{thm: main}]
    One direction is straightforward: an algebraic space $Y$ is smooth and proper if and only if $\Perf(Y)$ is smooth and proper (\cite[Theorem 11.4.0.3]{SAG}). Additionally, by GAGA (\cite[Corollary 7.5]{porta2016higher}, see also \cite[Example 5.3]{anmap}), we obtain the corresponding result for $\Perf(Y^{\mathrm{an}})$. The converse direction follows by combining the open immersion from \Cref{prop: X open into the analytification} with the algebraization result of \Cref{cor: algebraization of opens}.
\end{proof}
\subsection{Formal setting}
In this section, we consider formal schemes over $\Spf(\m{O}_{k})$, where $\m{O}_{k}$ is a discrete valuation ring (DVR) with maximal ideal $\mathfrak{m}$ and fraction field $k$. We assume that $\m{O}_k$ is complete with respect to the $\mathfrak{m}$-adic topology.
All formal schemes are assumed to be separated, quasi-compact, topologically almost of finite presentation, and locally Noetherian. Under these hypotheses, we rely on \cite{SAG} and \cite{antonio2019derived} as the main references for the theory of formal schemes in the derived setting.

Given a formal scheme $\F{X}$ as above, we can consider its \cat of perfect complexes $\Perf(\F{X})$. This is an idempotent-complete \cat, and by the Grothendieck Existence Theorem (see \cite[Theorem 8.5.0.3]{SAG}), it is also $\m{O}_{k}$-linear (i.e., it belongs to $\mathrm{Cat}^{\Perf}_{\m{O}_{k}}$).

\begin{definition}
    Let $X$ be a (derived) stack over $\Spec(\m{O}_{k})$. Its \emph{formalization} (also called formal completion) is defined as the following pullback: 
    \begin{equation*}
    \begin{tikzcd}
        {\widehat{X}} && X \\
        \\
        {\Spf(\m{O}_{k})} && {\Spec(\m{O}_{k})}
        \arrow[from=1-1, to=1-3]
        \arrow[from=1-1, to=3-1]
        \arrow["\lrcorner"{anchor=center, pos=0.125}, draw=none, from=1-1, to=3-3]
        \arrow[from=1-3, to=3-3]
        \arrow[from=3-1, to=3-3]
    \end{tikzcd}
\end{equation*}  
\end{definition}

In this section, we will prove an analogue of \Cref{thm: main} in the setting of formal schemes.

\begin{definition}
    Given a formal scheme $\F{X}$ over $\Spf(\m{O}_{k})$, we say that $\F{X}$ is \emph{algebraizable} if there exists an algebraic space $X$ over $\Spec(\m{O}_{k})$ and an equivalence $\F{X} \simeq \widehat{X}$. 
    In this case, we say that $\widehat{X}$ is the formal completion of $X$.
\end{definition}

The aim of this section is to prove the following theorem. 

\begin{theorem}\label{thm: main formal setting}
    Let $\F{X}$ be a smooth and proper formal scheme over $\Spf(\m{O}_{k})$. Then $\Perf(\F{X})$ is smooth and proper as an $\m{O}_k$-linear stable \cat if and only if $\F{X}$ is algebraizable. 
\end{theorem} 

In order to prove this theorem, we require some preliminary results. At the end of the section, we will also describe related results that do not rely on the smoothness hypothesis for $\F{X}$.\\

We begin by proving the following lemma, which extends \cite[Proposition 2.5]{liu} to the case of algebraic spaces (which we always assume to be separated and of finite type).

\begin{lemma}\label{lemma: extension of formal model}
    Let $X_{k}$ be an algebraic space over $k$. Then $X_{k}$ admits a model over $\m{O}_{k}$. Furthermore, if $X_k$ is proper, we can choose the model to be proper. 
\end{lemma}
\begin{proof}
We observe that, since $X_k$ is of finite type, we can extend it to a model $X$ over an open subset $V\subseteq \Spec(\m{O}_{k})$, such that $X$ is again separated and of finite type. Consider an atlas $U_{k}\to X_{k}$. By \cite[Proposition 2.5]{liu}, we can extend the scheme $U_{k}$ to a model $U$ that is faithfully flat over $\Spec(\m{O}_{k})$. By suitably shrinking $V$, we may assume that the restriction $U|_{V} \to X$ is a surjective \'etale map. Consider the following pushout diagram: 
    \begin{equation*}
    \begin{tikzcd}
	& {U|_{V}} & X \\
	& U & Y \\
	{} && {}
	\arrow[two heads, from=1-2, to=1-3]
	\arrow[hook, from=1-2, to=2-2]
	\arrow["j", from=1-3, to=2-3]
	\arrow["\alpha"', from=2-2, to=2-3]
	\arrow["\lrcorner"{anchor=center, pos=0.125, rotate=180}, draw=none, from=2-3, to=1-2]
\end{tikzcd}
    \end{equation*}
    By \cite[Lemma 3.1.6]{nagata}, we deduce that $j$ is injective, $\alpha$ is an atlas, and $Y$ is an algebraic space of finite type, faithfully flat over $\Spec(\m{O}_{k})$. By construction, $Y$ is a model of $X_{k}$, though it may not be proper. To resolve this, we invoke \cite[Theorem 1.2.1]{nagata} to find a proper algebraic space $\bar{Y}$ equipped with an open immersion $Y\to \bar{Y}$. The closure of $X$ inside $\bar{Y}$ then provides the required proper model, by \cite[Tag 0CMK]{stacks-project}. 
\end{proof}

We denote by $(-)^{\rig}$ the generic fiber functor for derived formal schemes and $\m{O}_{\F{X}}$-modules, as constructed in \cite{Ok-adicgeo}. In \emph{loc. cit.}, the authors establish the existence of a functor $(-)^{\rig}$ from the \cat of derived formal schemes to the \cat of $k$-analytic spaces, which, at the underived level, recovers the classical generic fiber functor of \cite{bosch2014lectures}. This functor naturally extends to categories of sheaves, leading to the following theorem:

\begin{theorem}\label{thm: rigidification of sheaves}
    Let $\F{X}$ be a formal scheme. The generic fiber functor induces a functor between the \cates of almost perfect complexes (still denoted by $(-)^{\rig}$)
    \begin{equation*}
        (-)^{\rig} \colon \mathrm{Coh}^{+}(\F{X}) \to \mathrm{Coh}^{+}(\F{X}^{\rig})
    \end{equation*}
    which preserves perfect complexes and is t-exact. Furthermore, for every $\m{F}\in \mathrm{Coh}^{b}(\F{X}^{\rig})$, there exists a formal model $\m{G} \in \mathrm{Coh}^{b}(\F{X})$ (i.e., $(\m{G})^{\rig} \simeq \m{F}$), and the category of formal models for $\m{F}$ is filtered.   
\end{theorem}
\begin{proof}
This is essentially the content of \cite[Appendix A.2.1]{Ok-adicgeo}. The proof demonstrates that if $\m{F}$ is bounded, we can choose $\m{G}$ to be bounded as well. The only remaining assertion is that the rigidification functor preserves perfect complexes. Working locally, we may assume $\F{X} = \Spf(A)$ for an $\m{O}_{k}$-algebra $A$. By \cite[Theorem 3.1]{derhom}, we can identify $\mathrm{Coh}^{+}((\Spf(A))^{\rig})$ with $\mathrm{Coh}^{+}(A\otimes_{\m{O}_{k}} k)$. Additionally, using \cite[Theorem 8.5.0.3]{SAG}, we can identify $\mathrm{Coh}^{+}(\Spf(A))$ with $\mathrm{Coh}^{+}(A)$. In this affine setting, the rigidification functor for sheaves coincides with the ordinary tensor product of modules: 
\[
-\otimes_{\m{O}_{k}}k \colon \Mod_{A} \to \Mod_{A\otimes_{\m{O}_{k}}k},
\]
which inherently preserves perfect complexes.
\end{proof}

\begin{lemma}\label{lem: smooth formal}
    Let $\F{X}$ be a formal scheme that is smooth over $\Spf(\m{O}_{k})$. Then $\Perf(\F{X})$ is equivalent to $\Cohb(\F{X})$.
\end{lemma}
\begin{proof}
Note that both \cates satisfy \'etale descent, so we may locally assume $\F{X} = \Spf(A)$. Furthermore, since $\F{X}$ is smooth, we may assume \'etale-locally that $A = \m{O}_k\langle T_{1},\dots,T_{n}\rangle$, implying that $\Spec(A)$ is regular over $\Spec(\m{O}_k)$ (see \cite[Proposition 7.3.23]{gabber2002ringtheorysixth}). In the locally Noetherian setting, the equivalence established in \cite[Proposition 8.4.1.2]{SAG} restricts to an equivalence of perfect objects and preserves t-structures (see \cite[Remark 8.4.2.6]{SAG}). In particular, this yields an equivalence between $\Cohb(\Spf(A))$ and $\Cohb(\Spec(A))$. Since $\Cohb(\Spec(A))$ is equivalent to $\Perf(\Spec(A))$, it follows that $\Perf(\Spf(A))$ is equivalent to $\Cohb(\Spf(A))$.
\end{proof}

\begin{corollary}\label{cor: smooth non archimedean}
    If $\F{X}$ is smooth over $\Spf(\m{O}_k)$, then $\Cohb(\F{X}^{\rig})$ is equivalent to $\Perf(\F{X}^{\rig})$.
\end{corollary}
\begin{proof}
    If $\m{F}$ is a bounded coherent complex over $\F{X}^{\rig}$, then by \Cref{thm: rigidification of sheaves} it admits a bounded formal model. This formal model is a perfect complex by \Cref{lem: smooth formal}. Consequently, its rigidification $\m{F}$ is also a perfect complex, again by \Cref{thm: rigidification of sheaves}. 
\end{proof}

\begin{definition}[{\cite[Definition 7.1.1.6]{SAG}}]
    Let $\m{C}$ be an $\m{O}_{k}$-linear \cat. We say that an object $C \in \m{C}$ is \emph{nilpotent relative to $\mathfrak{m}$} if, for every element $x \in \mathfrak{m}$, the relative tensor product 
    \[
    \m{O}_{k}[x^{-1}] \otimes_{\m{O}_{k}} C
    \]
    vanishes. We denote by $\m{C}_{nil}$ the full subcategory spanned by nilpotent objects.
\end{definition}

Informally, we are requiring that $\m{C}$ is supported set-theoretically at the closed point associated with the maximal ideal $\mathfrak{m}$. When $\m{C}$ is the \cat $\Cohb(\F{X})$ for a formal scheme $\F{X}$, our definition coincides with \cite[Definition 3.1]{derivedhilbert} (see \cite[Lemma 3.2]{derivedhilbert} and \cite[Example 7.1.1.2]{SAG}).

In the following lemma, we utilize one of the main theorems of \cite{thomason2007higher}. It states that there is a fiber-cofiber sequence in $\catperf$ of the form:
\begin{equation}\label{eq: fiber-cofiber}
    \Perf_{\mathfrak{m}}(\m{O}_{k}) \to \Perf(\m{O}_{k}) \to \Perf(k).
\end{equation}
Here, $\Perf_{\mathfrak{m}}(\m{O}_{k})$ denotes the \cat of perfect complexes supported on the closed point associated with the maximal ideal $\mathfrak{m}$ (i.e., complexes that vanish on the open complement of the point).

\begin{lemma}\label{lem: cofiber in catperf}
    Let $\m{A}$ be an \cat in $\catperf_{\m{O}_{k}}$. Then the map $\Perf_{\mathfrak{m}}(\m{O}_{k})\otimes_{\m{O}_{k}}\m{A} \to \m{A}$, obtained by tensoring the first map of \eqref{eq: fiber-cofiber} with $\m{A}$, is fully faithful, and its essential image coincides with $\m{A}_{nil}$. Furthermore, there is a cofiber sequence in $\catperf$ of the form:
    \begin{equation*}
        \Perf_{\mathfrak{m}}(\m{O}_{k})\otimes_{\m{O}_{k}}\m{A} \to \m{A} \to \Perf(k)\otimes_{\m{O}_{k}}\m{A}.
    \end{equation*}
\end{lemma}
\begin{proof}
By \cite[Proposition 7.1.5.3]{SAG}, we can identify $\Perf_{\mathfrak{m}}(\m{O}_{k})$ with $\Perf(\m{O}_{k})_{nil}$. The first statement for the ind-completion of our \cates follows from \cite[Corollary 7.1.2.11]{SAG}. By \cite[Lemma 3.10]{binda2021gaga}, we know that this map preserves compact objects, thus allowing us to extend the result to idempotent-complete \cates. To obtain the cofiber sequence, we simply tensor \eqref{eq: fiber-cofiber} by $-\otimes_{\m{O}_{k}}\m{A}$.    
\end{proof}

\begin{proposition}\label{prop: comparison perfect complexes}
    Let $\F{X}$ be a smooth and proper formal scheme over $\Spf(\m{O}_k)$. Then there is an equivalence of $k$-linear $\infty$-categories between $\Perf(\F{X})\otimes_{\m{O}_{k}}\Perf(k)$ and $\Perf(\F{X}^{\rig})$. 
\end{proposition}
\begin{proof}
    By the preceding lemma, we have a cofiber sequence:
    \begin{equation*}
        \Perf_{\mathfrak{m}}(\m{O}_{k})\otimes_{\m{O}_{k}}\Perf(\F{X}) \to \Perf(\F{X}) \to \Perf(k)\otimes_{\m{O}_{k}}\Perf(\F{X}).
    \end{equation*}
    We can identify $\Perf_{\mathfrak{m}}(\m{O}_{k})\otimes_{\m{O}_{k}}\Perf(\F{X})$ with $\Perf(\F{X})_{nil}$. \cite[Definition 3.5]{derivedhilbert} introduces the exact sequence: 
    \[
     \Cohb(\F{X})_{nil} \to \Cohb(\F{X}) \to \Cohb(\F{X})_\mathrm{loc}.
    \]
    \cite[Corollary 3.8]{derivedhilbert} proves that the right-hand side is equivalent to $\Cohb(\F{X}^{\rig})$ (note that all functors here are t-exact). By invoking \Cref{lem: smooth formal} and \Cref{cor: smooth non archimedean}, we obtain the equivalence stated in the proposition.
\end{proof}

\begin{proposition}\label{prop: algebrization}
	Let $\F{X}$ be a proper formal scheme such that $\F{X}^{\rig}$ is smooth and algebraizable. Then $\F{X}$ is algebraizable.
\end{proposition}
\begin{proof}
    By \Cref{thm: main}, there exists a smooth and proper algebraic space $Y$ over $\Spec(k)$ such that $Y^{\mathrm{an}} \simeq \F{X}^{\rig}$. Using \Cref{lemma: extension of formal model}, we can find a proper algebraic space $W$ over $\m{O}_{k}$ such that $\widehat{W}$ is a formal model for $\F{X}^{\rig}$. We now possess two formal models for $\F{X}^{\rig}$. By \cite[Theorem 4.1]{Bosch1993}, there exists another formal scheme $\tilde{Y}$ fitting into a zig-zag of formal schemes: 
    \begin{equation}
        \begin{tikzcd}
	&& \tilde{Y} \\
	& \F{X} && \widehat{W} \\
	{} && {}
	\arrow[from=1-3, to=2-2]
	\arrow[from=1-3, to=2-4]
    \end{tikzcd}
    \end{equation}
    where we may assume both maps are formal blow-ups. Since the map $\tilde{Y}\to \widehat{W}$ is a formal blow-up and $W$ is algebraic, we deduce by \cite[Theorem 3.1]{fmp} that $\tilde{Y}$ is also algebraic. Finally, since $\tilde{Y}$ is algebraic, applying \cite[Theorem 3.2]{fmp} allows us to conclude that $\F{X}$ must be algebraic as well.
\end{proof}

\begin{proof}[Proof of \Cref{thm: main formal setting}]
Assuming that $\Perf(\F{X})$ is smooth and proper, \Cref{prop: comparison perfect complexes} implies that $\Perf(\F{X}^{\rig})$ is smooth and proper over $k$, as it is the base change of a smooth and proper $\infty$-category. Using \Cref{thm: main}, we see that $\F{X}^{\rig}$ is algebraizable. \Cref{prop: algebrization} then guarantees that $\F{X}$ is algebraizable. The converse direction follows immediately from \cite[Theorem 8.5.0.3]{SAG}.
\end{proof}

\begin{remark}
	Note that to prove that $\F{X}$ is algebraizable, it is sufficient to assume that $\F{X}$ is proper, $\F{X}^{\rig}$ is smooth, and $\Perf(\F{X})$ is smooth and proper. Under these relaxed assumptions, the proof proceeds in exactly the same manner.  
\end{remark}

Using these same ideas, we can also prove the following result:

\begin{theorem}\label{thm: algebraization via coherent categories}
	Let $\F{X}$ be a proper formal scheme such that $\Cohb(\F{X})$ is smooth as a \cat in $\catperf_{\m{O}_{k}}$, and $\F{X}^{\rig}$ is smooth. Then $\F{X}$ is algebraizable. 
\end{theorem}
\begin{proof}
	Applying \Cref{lem: cofiber in catperf} with $\m{A}=\Cohb(\F{X})$, we obtain the cofiber sequence: 
    \begin{equation*}
        \Perf_{\mathfrak{m}}(\m{O}_{k})\otimes_{\m{O}_{k}}\Cohb(\F{X}) \to \Cohb(\F{X}) \to \Perf(k)\otimes_{\m{O}_{k}}\Cohb(\F{X}),
    \end{equation*}
    along with an equivalence between $\Perf_{\mathfrak{m}}(\m{O}_{k})\otimes_{\m{O}_{k}}\Cohb(\F{X})$ and $\Cohb(\F{X})_{nil}$. Exactly as in \Cref{prop: comparison perfect complexes}, we obtain an equivalence between $\Cohb(\F{X}^{\rig})$ and $\Cohb(\F{X})\otimes_{\Perf(\m{O}_{k})}\Perf(k)$. Since $\F{X}^{\rig}$ is smooth, \Cref{cor: smooth non archimedean} allows us to identify $\Cohb(\F{X}^{\rig})$ with $\Perf(\F{X}^{\rig})$. In particular, $\Cohb(\F{X}^{\rig})$ is smooth, and it is also proper by \Cref{cor: non commutative properness}. Thus, by \Cref{thm: main} and \Cref{prop: algebrization}, we conclude that $\F{X}$ is algebraizable.
\end{proof}

\begin{remark}
In particular, this theorem indicates that for a general proper formal scheme $\F{X}$, we do not expect $\Cohb(\F{X})$ to be smooth. For instance, consider a proper, non-algebraizable formal scheme $\F{X}$ whose rigid generic fiber is smooth (but not algebraizable). The preceding theorem forces us to conclude that $\Cohb(\F{X})$ cannot be smooth.

By contrast, for a separated scheme $X$ over a perfect field, Lunts demonstrated in \cite{lunts2010categorical} that $\Cohb(X)$ is always smooth. More generally, Efimov proved in \cite{efimov} that under mild conditions on $X$, the \cat $\Cohb(X)$ is of finite type. Because being of finite type implies smoothness, the above proposition highlights a sharp contrast: neither of these statements can be generalized unconditionally to the setting of formal schemes.
\end{remark}

\printbibliography

\Addresses

\end{document}